\magnification=\magstep1
\input amstex
\documentstyle{amsppt}
\hsize = 6.5 truein
\vsize = 9 truein

\UseAMSsymbols

\define \J {\Bbb J}
\define \rL {\Bbb L}
\define \F {\Bbb F}
\define \K {\Bbb K}
\define \rO {\Bbb O}
\define \N {\Bbb N}
\define \AN {\Bbb N^*}
\define \rP {\Bbb P}
\define \CP {\Cal P}
\define \CET {\Cal C (ET)}
\define \s {\Bbb S}
\define \zs {\Bbb S_0}
\define \CRET {\Cal C (RET)}
\define \Z {\Bbb Z}
\define \Zi {\Bbb Z[i]}
\define \Q {\Bbb Q}
\define \R {\Bbb R}
\define \C {\Bbb C}
\define \oN {\overline {\Bbb N}}
\define \CPET {\Cal C (PET)}
\define \CePET {\Cal C_e (PET)}
\define \oR {\overline {\Bbb R}}
\define \RP {\R^{+*}}
\define \cal {\text {cal}}
\define \th	{\theta}
\define \cp {\text {cap}}
\define \gri {\text {gri}}
\define \lgr {\text {lgr}}
\define \grs {\text {grs}}
\define \ugr {\text {ugr}}
\define \gr {\text {gr}}
\define \e	{\epsilon}
\define \de	{\delta}
\define \p {\Bbb P}
\define \Gi {\Bbb Z[i]}

\define \ee {\'e}
\define \ea {\`e}
\define \oo {\^o}

\define \CC {\Cal C}
\define \CDD {\Cal D}
\define \CF {\Cal F}
\define \CG {\Cal G}
\define \CH {\Cal H}
\define \CI {\Cal I}
\define \CL {\Cal L}
\define \CN {\Cal N}
\define \CT {\Cal T}
\define \CU {\Cal U}
\define \CV {\Cal V}
\define \CW {\Cal W}
\define \CA {\Cal A}
\define \CB {\Cal B}
\define \CE {\Cal E}
\define \CR {\Cal R}
\define \CS {\Cal S}
\define \CM {\Cal M}
\define \CO {\Cal O}
\define \CX {\Cal X}
\define \CZ {\Cal Z}

\define \Card {\text{Card}}

\define \bS {\bf S \rm}

\define \noo {$\text{n}^o$}

\define \T- {\overset {-1}\to T}
\define \G- {\overset {-1}\to G}
\define \GA- {\overset {-1}\to \Gamma}
\define \RA- {\overset {-1}\to R}
\define \PA- {\overset {-1}\to P}
\define \ig {g^{-1}}

\define \Ž {\'e}
\define \ {\`e}
\define \ˆ {\`a}
\define \ {\`u}

\define \™ {\^o}
\define \ {\^e}
\define \" {\^\i}
\define \ž {\^u}
\define \‰ {\^a}

\define \rd {\roman d}
\define \rD {\roman D}
\define \rad {\roman {ad}}
\define \rp {\roman p}
\define \rs {\roman s}
\define \Der {\roman {Der}}

\define \su {\subheading}
\define \noi {\noindent}

\UseAMSsymbols
\topmatter
\title Entrelacement d'alg\`ebres de Lie \endtitle

\author Barben-Jean COFFI-NKETSIA
\footnote"*"{\`A \sl Coffi\rm, d\Ž c\Ž d\Ž
\ pendant la pr\Ž paration et avant l'ach\ vement de ce
texte, un amical adieu. \sl Haddad\rm} et Labib
HADDAD\endauthor

\leftheadtext{B.-J. COFFI-NKETSIA et L. HADDAD}

\address {120 rue de Charonne, 75011 Paris, France;
e-mail: labib.haddad\@wanadoo.fr} \endaddress
 
\abstract\nofrills{\smc English abtract.} \bf Wreath
product for Lie algebras\rm. Full details are given for the
definition and construction of the wreath product $W(A,B)$
of two Lie algebras $A$ and
$B$, in the hope that it can lead to the definition of a
suitable Lie group to be the wreath product of two given Lie
groups. In the process, quite a few new notions
are needed, and introduced. Such are, for example : Formal
series with variables in a vector space and coefficients in
some other vector space. Derivation of a formal series
relative to another formal series. The Lie algebra of a
vector space. Formal actions of Lie algebras over vector
spaces. The basic formal action of a Lie algebra over
itself (as a formal version of the analytic aspect of the
infinitesimal operation law of a Lie groupuscule). More
generally, relative to any given formal action
$\rd$ of the Lie algebra $B$ over a vector space $Y$, the
wreath product $W(A,B;\rd)$ is defined. When $\rd$ is the
basic action of $B$ over itself, the special case
$W(A,B)$ is recovered. Main features are : A description of
the triangular actions of the wreath products
$W(A,B;\rd)$ over product vector spaces $X \times Y$. A
Kaloujnine-Krasner type theorem : In essence, it says that
all Lie extensions $C$ of the Lie algebra $B$ by the Lie
algebra $A$ are, indeed, subalgebras of their wreath product
$W(A,B)$. A moderately detailed english summary of the paper
can be found in [7, p. 9-15].
\endabstract
\endtopmatter
\document
\refstyle{A}
\widestnumber\key{ABCD}
\su{0} Introduction.
\su{1} S\Ž ries formelles.
\su{2} Sym\Ž trisation.
\su{3} D\Ž rivation suivant une s\Ž rie formelle.
\su{4} L'alg\ bre de Lie $S(X)$ d'un espace vectoriel $X$.
\su{5} Action d'une alg\ bre de Lie sur un espace vectoriel.
\su{6} L'exemple originel.
\su{7} Produits d'entrelacements.
\su{8} Action triangulaire.
\su{9} Action fondamentale d'une alg\ bre de
Lie sur elle-m\ me.
\su{10} Produit d'entrelacement de deux alg\ bres de Lie.
\su{11} Repr\Ž sentation des extensions dans le produit
d'entrelacement.
\newpage

\
\head 0. Introduction\endhead

\

Le produit d'entrelacement (\lq\lq wreath
product" en anglais) de deux groupes (abstraits)
quelconques,
$A$ et
$B$, est une notion devenue classique : voir, par exemple, 
Neumann [6]. Voici, succintement, comment on l'introduit. 

\

 Consid\Ž r\Ž \ comme
produit de copies du groupe $A$, 
l'ensemble $A^B$ de toutes les applications de $B$ dans
$A$ est, lui-m\ me, un groupe. Soit
$\sigma : B
\times A^B \to A^B$ l'application $(b,f) \mapsto
\sigma_b(f)$ o\ \ $\sigma_b(f) \in A^B$ est l'application
de $B$ dans $A$ d\Ž finie par $(\sigma_b(f))(x) = f(x \
b)$ pour tout
$x
\in B$. Chacun des $\sigma_b$ est alors un automorphisme du
groupe $A^B$ et l'on v\Ž rifie que l'application $\sigma : B
\to \roman{Aut}(A^B)$ ainsi d\Ž finie est un homomorphisme
du groupe $B$ dans le groupe $\roman{Aut}(A^B)$ des
automorphismes du groupe
$A^B$. 

\

Le produit d'entrelacement de $B$ par $A$ est alors
d\Ž fini comme \Ž tant le produit semi-direct  $W(A,B) =
A^B \underset\sigma\to \times B$ du groupe $B$ par le groupe
$A^B$ suivant l'homomorphisme $\sigma$. Plus pr\Ž cis\Ž ment,
pour deux
\Ž l\Ž ments
$(f,b)$ et
$(g,c)$  de
$W(A,B)$, on a
$$(f,b) \ (g,c) = (f \ \sigma_b(g),b \ c).$$
On introduit une op\Ž ration, dite
\sl triangulaire\rm,
 du groupe
$W(A,B)$ sur le produit cart\Ž sien
$A \times B$, naturellement, en posant $(x,y).(f,b) = (x \
f(y),y\ b)$ pour chaque
op\Ž rateur $(f,b) \in A^B
\times B$ et chaque point $(x,y) \in A \times B$. Ce
faisant, on fait ainsi op\Ž rer ce groupe  \sl \ˆ
\ droite\rm.

\

L'une des propri\Ž t\Ž s du produit d'entrelacement est
connue sous le nom de th\Ž or\ me de Kaloujnine-Krasner :
tout groupe $C$ extension de
$B$ par
$A$ est un sous-groupe du produit d'entrelacement $W(A,B)$.

\

M. Krasner a pos\Ž \ le probl\ me suivant au premier des
deux auteurs du texte pr\Ž sent : \bf Comment d\Ž finir un
produit d'entrelacement dans le cas des groupes de Lie?\rm

\

Autrement dit, lorsque $A$ et $B$ sont des groupes de Lie,
comment peut-on distinguer un sous-groupe particulier de
$W(A,B)$, assez \sl convenable \rm
pour en faire un groupe de Lie, et de mani\ re
\sl raisonnable\rm.

\

Le probl\ me semblait assez ardu, difficile et hors de
port\Ž e. A notre connaissance, il n'a d'ailleurs toujours
pas
\Ž t\Ž \ r\Ž solu. Aussi, nous sommes-nous attel\Ž , tout
d'abord, au probl\ me qui nous a paru plus simple, celui
du produit d'entrelacement d'alg\ bres de Lie. On passait
ainsi du domaine \sl analytique et g\Ž om\Ž trique\rm, tr\ s
riche, au domaine
\sl alg\Ž brique \rm plus accessible. Il n'\Ž tait pas trop
hasardeux de penser que l'alg\ bre de Lie d'un \sl \Ž
ventuel \rm groupe de Lie qui serait un produit
d'entrelacement de deux groupes de Lie donn\Ž s serait le
produit d'entrelacement des alg\ bres de Lie de ces deux
groupes ou, pour le moins, une sous-alg\ bre \sl
convenable \rm de ce produit d'entrelacement.

\

Partis en qu\ te du \sl v\Ž ritable \rm produit
d'entrelacement de deux alg\ bres de Lie, nous avons suivi
un chemin sinueux dont il serait trop long, ici, de d\Ž crire
en d\Ž tail les \Ž tapes. Qu'il suffise de dire ceci, pour
le moment.

\

D'un c\™ t\Ž , la route montait vers un hypoth\Ž tique
produit d'entrelacement de deux groupes de Lie, entrevu, \ˆ \
la mani\ re des groupes abstraits, de l'autre, elle
redescendait vers les alg\ bres de Lie correspondantes dont
on essayait de deviner les structures et l'enlacement...

\

Nous avons rencontr\Ž , \ˆ \ mi-chemin, entre l'analytique
et l'alg\Ž brique, la notion de loi d'op\Ž ration infinit\Ž
simale et celle de groupuscule de Lie. Apr\ s de multiples 
et infructueux essais, elles nous ont permis de d\Ž gager la
voie et d'ouvrir le passage. Le Trait\Ž \ de Bourbaki, et
ses tr\ s pr\Ž cieux exercices, nous ont grandement
facilit\Ž \ la t\‰ che.

\

On passerait directement du cas des groupes abstraits \ˆ \
celui des alg\ bres de Lie, tout simplement, en rempla\c
cant le groupe $A^B$ par une alg\ bre de Lie $A[[B]]$,
convenable, et en faisant jouer le r\™ le du groupe
$\roman{Aut}(A^B)$ des automorphismes du groupe $A^B$ par
l'alg\ bre de Lie des d\Ž rivations de $A[[B]]$.

\

Pour cela, il suffisait de trouver les \sl bonnes \rm d\Ž
finitions pour $A[[B]]$ et sa structure d'alg\ bre de Lie.
Ce que nous f\^\i mes, en choisissant pour $A[[B]]$ le statut
d'une alg\ bre de Lie de s\Ž ries formelles.

\ 

Dans le texte qui suit, nous donnons toutes les d\Ž
finitions et les d\Ž tails n\Ž cessaires pour la construction
du produit d'entrelacement de deux alg\ bres de Lie, sur un
corps commutatif, de caract\Ž ristique nulle, quelconque.

\

Dans un texte ant\Ž rieur, \lq\lq Produit d'entrelacement et
action triangulaire d'alg\`ebres de
Lie" [7], nous avons donn\Ž \
un r\Ž sum\Ž
\ assez complet de ce qui va suivre. Pour la concordance de
nos deux textes, nous utilisons, de nouveau, les num\Ž ros de
renvoi sous la forme $<n>$.

\

Les d\Ž tails sont nombreux avant d'arriver au
produit d'entrelacement de deux alg\ bres de Lie $A$ et
$B$ sous la forme $W(A,B) = A[[B]] \times B$. 

\

Parmi ces d\Ž tails, on
signalera l'introduction de quelques notions nouvelles. En
parti\-culier, l'ensemble $X[[Y]]$ des s\Ž ries formelles, \ˆ
\ variables dans un espace vectoriel $Y$, donn\Ž , et
coefficients dans un autre espace vectoriel $X$, et dont
$A[[B]]$ est un cas particulier.  La d\Ž rivation d'une s\Ž
rie formelle suivant une autre. L'alg\ bre de Lie
$S(X) = X[[X]]$ d'un espace vectoriel $X$ quelconque. Les
actions des alg\ bres de Lie sur les espaces
vectoriels. L'action fondamentale d'une alg\ bre de Lie sur
elle-m\ me (version formelle de l'aspect analytique de la
loi d'op\Ž ration infinit\Ž simale d'un groupuscule de Lie).

\

On trouvera \Ž galement, ici, la d\Ž finition d'un
produit d'entrelacement plus g\Ž n\Ž ral, $W(A,B;\rd) =
A[[Y]]
\times B$,  relatif
\ˆ
\ une action formelle quelconque $\rd$ de l'alg\ bre de Lie
$B$ sur un espace vectoriel donn\Ž \ $Y$. Lorsque $\rd$ est
l'action fondamentale de $B$ sur elle-m\ me, on retrouve le
cas particulier du produit $W(A,B)$.

\

Signalons les deux points saillants suivants. Une
description des actions
\sl triangulaires \rm des produits d'entrelacement
$W(A,B;\rd)$ sur les espaces vectoriels produits $X \times
Y$. Ainsi qu'un th\Ž or\ me du type Kaloujnine-Krasner :
toute alg\ bre de Lie $C$ extension de $B$ par $A$ se
plonge dans l'alg\ bre de Lie $W(A,B)$. Ce plongement, loin
d'\ tre banal, n'est pas unique. On verra, au dernier
paragraphe, la formule remarquable de ces plongements.

\

Un dernier mot encore au sujet du probl\ me pos\Ž \
initialement par M. Krasner. Dans le cas r\Ž el ou
complexe, m\ me lorsque $A$ et $B$ sont, toutes deux,  de
dimension finie, il sera sans doute malais\Ž \ de remonter
directement du produit d'entrelacement de ces deux alg\ bres
de Lie
\ˆ \ un groupe de Lie $G$ dont l'alg\ bre de Lie serait
$L(G) = W(A,B)$. En effet, la dimension de $W(A,B)$ est
toujours infinie (sauf si $A$ ou $B$ est triviale). On
pourra probablement remonter jusqu'\ˆ \ un groupuscule de
Lie. Mais, au-del\ˆ , il faudrait prendre pour $G$ un
groupe de Lie banachique, port\Ž \ par une vari\Ž t\Ž \
analytique model\Ž e sur un espace de Banach (voire de
Hilbert) de dimension infinie, par exemple, ou s\Ž lectionner
une sous-alg\ bre de dimension finie, \sl convenable\rm, de
l'alg\ bre de Lie
$W(A,B)$. Tout cela demeure en dehors du champ de ce modeste
essai.

\

\

\

\bf Tout au long de ce texte\rm, $K$
d\Ž signe un corps commutatif
\bf infini \rm quelconque; il sera de caract\Ž ristique
nulle \ˆ \ partir du paragraphe \bf 9\rm. De m\ me,
$E,F,X,Y,$ d\Ž signeront des espaces vectoriels sur $K$, et
$A,B,C,$ d\Ž signerons des $K$-alg\ bres de Lie. Par
$m,n,r,$ on d\Ž signera des entiers naturels quelconques,
autrement dit des entiers $\geq 0$. Pour chaque
$m$, on d\Ž signera par $L_m(E;F)$ l'ensemble de toutes les
applications $m$-lin\Ž aires de $E$ \ˆ \ valeurs dans $F$.

\

\bf Tous les espaces vectoriels ainsi que toutes les alg\
bres de Lie seront suppos\Ž s avoir $K$ comme corps des
scalaires\rm.

\

\

On a essay\Ž \ d'\ tre le plus complet possible, pour
facilter la lecture du texte, en donnant tous les d\Ž tails
utiles \ˆ
\ sa compr\Ž hension (frisant parfois, sans doute,  la
surcharge).

\

\

\

\newpage

\head 1. S\Ž ries formelles\endhead

\

\eightpoint Bourbaki, dans un court \smc APPENDICE \rm [1, p.
88-89], expose les notions de polyn\™ mes-continus et de s\Ž
ries formelles pour les espaces vectoriels polynorm\Ž s, en
passant par les applications \bf multilin\Ž aires
continues\rm. En imitant cette d\Ž marche, nous l'adaptons
ci-dessous au cas
purement alg\Ž brique des espaces vectoriels
quelconques (et la d\Ž veloppons pour nos besoins).\tenpoint

\

\su{1.1. Polyn\™ mes homog\ nes} On dira qu'une application
$f : E \to F$ est un \bf polyn\™ me homog\ ne \rm de degr\Ž
\ $m$, \ˆ \ variables dans $E$ et coefficients dans $F$,
lorsqu'il existe au moins une application $m$-\bf lin\Ž aire
\rm 
$u : E^m \to F$ telle que
$$f(x) = u(x,\dots,x) \ \ \text{pour tout} \ \ x \in E.$$
On dira alors que $f$ est \bf d\Ž termin\Ž \ \rm par $u$.

\

On d\Ž signera par $F[E]_m$ l'ensemble de ces polyn\™ mes
homog\ nes de degr\Ž \ $m$. C'est naturellement un espace
vectoriel sur $K$.

\

\su{1.2. Remarques}

\

\noi $1^o$. Le corps $K$ \Ž tant infini, la somme 
$$\sum_m F[E]_m,$$
consid\Ž r\Ž e dans le $K$-espace vectoriel $F^E$ de toutes
les applications de $E$ dans $F$, est une \bf somme
directe\rm.

\

\noi \bf En effet\rm, soit $f = f_0 + f_1 + \dots + f_m +
\dots + f_r$ o\ \ $f_m \in F[E]_m$ pour chaque $m$. On
suppose que $f(x) = 0$ pour tout $x \in E$. Alors, pour tous
$\lambda \in K$ et $x \in E$, on aura
$$f(\lambda x) = f_0(x) + \lambda f_1(x) + \dots + \lambda^r
f_r(x) = 0$$
et, $K$ \Ž tant infini, on aura donc
$$f_0(x) = f_1(x) = \dots = f_r(x) = 0 \ \ \text{pour tout}
\ \ x \in E.\qed$$

\

\noi $2^o$. \bf Contre-exemple\rm. Bien entendu, sur le
corps premier $\F_p$, les deux polyn\™ mes
$$x_1^p + x_2^p + \dots + x_n^p \ \ \text{de degr\Ž \ } \ 
p$$ 
et
$$x_1 + x_2 + \dots + x_n \ \ \text{de degr\Ž \ } \  1$$
sont \Ž gaux !

\

\su{1.3. Polyn\™ mes et s\Ž ries formelles}

\

On posera 
$$F[E] = \bigoplus_m F[E]_m \ , \ F[[E]] = \prod_m
F[E]_m.$$
On appelle alors \bf polyn\™ me \rm (resp. \bf s\Ž rie
formelle\rm) \ˆ \ variables dans $E$ et coefficients dans
$F$ tout \Ž l\Ž ment de $F[E]$ (resp. de $F[[E]]$).

\

Ainsi, une s\Ž rie formelle $f \in F[[E]]$ est de la forme
$$f = f_0 + f_1 + \dots + f_m + \cdots$$
o\ \ chaque $f_m$, appel\Ž \ \bf composante homog\ ne \rm
de degr\Ž \ $m$ de $f$, appartient \ˆ \ $F[E]_m$.

\

\su{1.3.1. Exemples primaires} Que sont les polyn\™ mes
homog\ nes et les s\Ž ries formelles dans le cas o\ \ les
deux espaces vectoriels $E$ et $F$ sont de dimensions finies?

\

Soient $E =
K^r$ et $F = K^s$.

\

\noi Lorsque $s = 1$, l'ensemble des polyn\™ mes
$F[E]$ s'identife \ˆ \ $K[x_1,\dots,x_r]$, ensemble des
polyn\™ mes (classiques) \ˆ \ $r$ variables et \ˆ \
coefficients dans le corps $K$. De m\ me, $F[[E]]$
s'identifie alors 
\ˆ
\
$K[[x_1,\dots,x_r]]$, ensemble des s\Ž ries formelles
(classiques) aux $r$ variables $x_1,\dots,x_r,$ \ˆ \
coefficients dans le corps $K$. Ainsi, pour $s$ quelconque,
chaque
$f \in F[[E]]$ s'identifie \ˆ \ un $s$-uple
$(f_1,\dots,f_s)$ de s\Ž ries formelles $f_i \in
K[[x_1,\dots,x_r]]$.

\

\noi Autrement dit, $F[[E]]$ s'identifie,
canoniquement, au produit de $s$ copies de
$K[[x_1,\dots,x_r]]$.
$$(K^s)[[K^r]] = (K[[x_1,\dots,x_r]])^s = \prod_{i=1}^s
K_i[[x_1,\dots,x_r]].$$

\

\su{1.3.2. Sommabilit\Ž }

\

On commence par observer ceci. Soit $(f_i)_{i \in I}$ une
famille de s\Ž ries formelles 
$$f_i \in F[[E]] \ \ , \ \ f_i = \sum_m f_{i,m} \ \ , \ \
f_{i,m}
\in F[E]_m.$$
Lorsque $I$, l'ensemble des indices, est fini, la somme
$\sum_{i \in I} f_i$ est bien d\Ž finie. Plus g\Ž n\Ž
ralement, pour $I$ quelconque, fini ou infini, soit
$$I(m) = \{ i \in I : f_{i,m} \neq 0 \}.$$
Si tous les $I(m)$ sont finis, on pose $g_m = \sum_{i \in
I(m)} f_{i,m}$, de sorte que $g_m \in F[E]_m$. On dit alors
que la famille est
\bf sommable
\rm et que sa somme est
$$\sum_{i
\in I}  f_i = g_0 + g_1 + \dots + g_m + \cdots\cdot$$

\

\su{1.3.3. Question de m\Ž thode} Dans la
suite, pour \Ž tablir certains des nombreux r\Ž sultats sur
les s\Ž ries formelles, on utilisera souvent la m\ me m\Ž
thode : elle consiste \ˆ \ \Ž tablir ce r\Ž sultat, d'abord
pour les polyn\™ mes homog\ nes (en utilisant les
applications multilin\Ž aires qui les d\Ž terminent) puis,
en l'appliquant
\ˆ
\ leurs composantes homog\ nes, \ˆ
\ l'\Ž tendre aux s\Ž ries formelles elles-m\ mes.

\

Dans la plupart des cas, on a donn\Ž \ tous les d\Ž
tails des calculs, assez arides et fastidieux, pour le
confort du lecteur, ainsi que pour notre propre tranquillit\Ž
 \dots

\

\su{1.4. S\Ž ries formelles doubles}

\

Comme dans le cas classique, on peut introduire la notion
plus g\Ž n\Ž rale de s\Ž rie formelle \bf double \rm \ˆ \
variables dans deux espaces vectoriels $X$ et $Y$, et \ˆ \
coefficients dans un troisi\ me espace $F$.

\

Voici comment on proc\ de.

\

On d\Ž signe par $L_{(m,n)}(X,Y;F)$ l'ensemble de toutes les
applications $u : X^m \times Y^n \to F$ qui sont
$(m+n)$-multilin\Ž aires, autrement dit, telles que
$u(x_1,\dots,x_m,y_1,\dots,y_n)$ soit lin\Ž aire par rapport
\ˆ \ chacune des variables $x_1,\dots,x_m,y_1,\dots,y_n$.

\ 

On dira qu'une application $f : X \times Y \to F$ est un
polyn\™ me homog\ ne de \bf bidegr\Ž \ \rm $(m,n)$, \ˆ \
variables dans $(X,Y)$ et coefficients dans $F$, lorsqu'il
existe au moins un \Ž l\Ž ment $u \in L_{(m,n)}(X,Y;F)$ tel
que
$$f(x,y) = u(x,\dots,x,y\dots,y) \ \ \text{pour tous} \ \ x
\in X \ , \ y \in Y.$$
On d\Ž signera par $F[X,Y]_{(m,n)}$ l'ensemble de ces
polyn\™ mes homog\ nes de bidegr\Ž \ $(m,n)$. C'est
naturellement un espace vectoriel sur $K$.

\

On posera $F[[X,Y]] = \dsize\prod_{(m,n)} F[X,Y]_{(m,n)}$.

\

On appellera s\Ž rie formelle double, \ˆ \ variables dans
$(X,Y)$ et coefficients dans $F$, tout \Ž l\Ž ment de
$F[[X,Y]]$.

\

\su{1.5. Identifications naturelles et canoniques}

\

On consid\ re $T = F[[X]]$ et on \Ž crit $T_m$ pour
$F[X]_m$.

\

On fera les identifications suivantes :
$$L_{(m,n)}(X,Y;F) \ \ \text{est identifi\Ž \ avec} \ \
L_n(Y;L_m(X;F)).\tag"\bf 1."$$
$$F[X,Y]_{(m,n)} \ \ \text{est identifi\Ž \ avec} \ \
T_m[Y]_n.\tag"\bf 2."$$
$$F[[X,Y]] \ \ \text{est identifi\Ž \ avec} \ \
T[[Y]].\tag"\bf 3."$$
$$F[[X \times Y] \ \ \text{est identifi\Ž \ avec} \ \
F[[X,Y]].\tag"\bf 4."$$
Ainsi
$$F[[X \times Y] \ \ \text{est identifi\Ž \ avec} \ \
T[[Y]].\tag"\bf 5."$$
Voici les d\Ž tails essentiels de ces identifications. Le
lecteur ne rencontrera pas de difficult\Ž \ \ˆ \ les compl\Ž
ter par des calculs, s'il le juge utile.

\

\su{1.5.1} Soit $u : X^m \times Y^n \to F$ une application
$(m+n)$-multilin\Ž aire. A chaque $y = (y_1,\dots,y_n) \in
Y^n$ correspond une application $u_y : X^m \to F$ d\Ž finie
par
$$u_y(x_1,\dots,x_m) = u(x_1,\dots,x_m,y_1,\dots,y_n).$$
Ainsi, $u_y \in L_m(X;F)$. Soit alors $\bar u$ l'application
de $Y^n$ dans $L_m(X;F)$ qui, \ˆ \ chaque $y \in Y^n$, fait
correspondre $u_y$, de sorte que $\bar u \in
L_n(Y;L_m(X;F))$. L'application $u \mapsto \bar u$ est
alors un
\bf isomorphisme \rm naturel, canonique, entre les deux
espaces vectoriels $L_{(m,n)}(X,Y;F)$ et $L_n(Y;L_m(X;F))$.

\

De cette premi\ re identification d\Ž coule simplement la
suivante.

\

\su{1.5.2} Soit $f \in F[X,Y]_{(m,n)}$ un polyn\™ me homog\
ne. Il est d\Ž termin\Ž \ par un \Ž l\Ž ment $u \in
L_{(m,n)}(X,Y;F)$. Le polyn\™ me $\bar f$ d\Ž termin\Ž \ par
$\bar u$ appartient donc \ˆ \ $T_m[Y]_n$. De plus $\bar f$
ne d\Ž pend que de $f$ et pas du choix de $u$ !

\

L'application $f \mapsto \bar f$ est un \bf isomorphisme \rm
entre les deux espaces vectoriels 
$$F[X,Y]_{(m,n)} \ \ \text{et} \ \ T_m[Y]_n.$$
Retenons les deux identit\Ž s suivantes :
$$\bar u(y_1,\dots,y_n)(x_1,\dots,x_m) = u(x_1,\dots,x_m,
y_1,\dots,y_n)$$
$$\bar f(y)(x) = f(x,y).$$

\

De cette seconde identification, on passe \ˆ \ la troisi\
me comme suit.

\

\su{1.5.3} Soit $f = \dsize\sum_{(m,n)} f_{(m,n)} \in
F[[X,Y]]$. Pour chaque $n$, soit $g_n = \dsize\sum_m \bar
f_{(m,n)}$, de sorte que $g_n \in T[Y]_n$.

\

On pose alors $\bar f = \sum_n g_n$ et on remarque que
cette notation est en accord avec celle qui a \Ž t\Ž \
adopt\Ž e ci-dessus dans le cas o\ \ $f$ est un polyn\™ me
homog\ ne. On a donc $\bar f \in T[[Y]]$, et l'application
$f \mapsto \bar f$ est un \bf isomorphisme \rm entre les
deux espaces vectoriels
$$F[[X,Y]] \ \ \text{et} \ \ T[[Y]].$$
Elle est, en effet, lin\Ž aire et \bf injective\rm. Il
suffit de se convaincre qu'elle est, \Ž galement, \bf
surjective\rm.

\

Pour cela, on se donne
$$g \in T[[Y]] \ , \ g = \sum_m g_m \ , \ g_m \in T[Y]_m.$$
Pour $y \in Y$, on a 
$$g_m(y) \in T = F[[X]] \ , \ g_m(y) =
\sum_n g_{m,n}(y) \ , \ g_{m,n}(y) \in F[X]_n.$$
Pour chaque $x \in X$, on a 
$$g_{m,n}(y)(x) \in F.$$
On pose alors
$$f_{m,n}(x,y) = g_{m,n}(y)(x).$$
On sait que l'on a
$$f_{m,n} \in F[X,Y]_{(m,n)} \ \ \text{et} \ \
g_{m,n} = \bar f_{m,n}.$$
Ainsi, en posant
$$f = \sum_{(m,n)} f_{(m,n)} \in F[[X,Y]],$$
on a
$$g = \bar f.$$

\

A nouveau, et formellement cette fois, on
peut
\Ž crire
$$\bar f(y)(x) = f(x,y).$$

\

\su{Remarque} Comme on l'a vu pour les polyn\™ mes simples
(ci-dessus \bf 1.2\rm), ici aussi  la somme 
$\dsize\sum_{(m,n)} F[X,Y]_{(m,n)}$ dans l'espace vectoriel
$ F^{X \times Y}$ est une \bf somme directe\rm.

\

\noi \bf En effet\rm,  soit $f = \sum_{m+n \leq r}
f_{(m,n)}$ o\ \ $f_{(m,n)} \in F[X,Y]_{(m,n)}$. On suppose
que $f(x,y) = 0$ pour tous $x \in X$ et $y \in Y$. On veut
montrer que $f_{(m,n)} = 0$ pour tout
$(m,n)$, autrement dit,
$f_{(m,n)}(x,y) = 0$ pour tous $x,y,m,n$. 

\

Pour chaque
$\lambda \in K \ , \ \mu \in K \ , \ x \in X \ , \ y \in Y$,
on a $f(\lambda x,\mu y) = 0$. Or,
$$f(\lambda x, \mu y) = \sum_{(m,n)} f_{(m,n)} (\lambda
x,\mu y) = \sum \lambda^m \mu^n f_{(m,n)}(x,y).$$
D'o\ \ la conclusion.\qed

\  

On en vient \ˆ \ la quatri\ me identification.

\

\su{1.5.4. Identification de $F[[X \times Y]$ avec
$F[[X,Y]]$} Elle r\Ž sultera de l'identification de
$$F[X \times Y]_r \ \ \text{avec} \ \ \bigoplus_{m+n = r}
F[X,Y]_{(m,n)} = F[X,Y]_r.$$
On pose $Z = X \times Y$ et on identifie $X$ au sous-espace
$X \times \{0\} \subset Z$, et $Y$ au sous-espace $\{0\}
\times Y \subset Z$. Soit $f \in F[X \times Y]_r$. 

\

On
consid\ re $\hat f(x,y) = f(x+y)$. 

\

\su{1.5.4.1} On montre que $\hat f \in F[X,Y]_r$.

\

\noi \bf En effet \rm : $f$ est d\Ž termin\Ž \ par un $u \in
L_r(Z;F)$ et, pour $z = x + y \ , \ x \in X \ , \ y \in Y,$
on a
$$f(z) = u(z,\dots,z) = u(x+y,\dots,x+y).$$
Or, 
$$u(z_1,\dots,z_r) = u(x_1 + y_1,\dots,x_r + y_r) \ \
\text{o\ \ } \ \ z_i = x_i + y_i.$$
Soit $R = \{1,\dots,r\}$. Pour chaque partie $M \subset R$,
on consid\ re le sous-espace $Z_M$ de $Z^r$ form\Ž \ des \Ž
l\Ž ments $z = (z_1,\dots,z_r)$ o\ \
$$z_i \in X \ \ \text{si} \ \ i \in M,$$
$$z_i \in Y \ \ \text{si} \ \ i \notin M.$$
La restriction de $u$ \ˆ \ $Z_M$ d\Ž finit ainsi
naturellement une application $(m+n)$-multilin\Ž aire
$$u_M : X^m \times Y^n \to F \ , \ m = \Card(M).$$
A son tour, $u_M$ d\Ž termine un polyn\™ me homog\ ne $f_M
\in F[X,Y]_{(m,n)}$,
$$f_M(x,y) = u_M(x,\dots,x,y,\dots,y).$$
Or
$$u(x + y,\dots,x + y) = \sum_M u_M(x,\dots,x,y\dots,y).$$
De sorte que $\hat f(x,y) = f(x+y) = \sum_M f_M(x,y),$
donc $\hat f \in F[X,Y]_r$.

\

Visiblement, l'application ainsi d\Ž finie de $F[X \times
Y]_r$ dans $F[X,Y]_r$, qui associe $\hat f$ \ˆ \ $f$, est \bf
lin\Ž aire\rm.

\

\su{1.5.4.2. Elle est aussi injective} Car si $\hat f \equiv
0$ alors $f(x+y) = \hat f(x,y) = 0$ pour tous $x,y,$
autrement dit,
$$f(z) = 0 \ \ \text{pour tout} \ \ z \in Z.$$
Donc $f = 0$.

\  

\su{1.5.4.3. Cette application est \Ž galement surjective}

\

\noi \bf En effet\rm, soit $f \in F[X,Y]_{(m,n)}$ d\Ž
termin\Ž \ par un $u \in L_{(m,n)}(X,Y;F)$. On d\Ž finit
alors $v : Z^r \to F$ par
$$v(z_1,\dots,z_r) = u(x_1,\dots,x_m,y_{m+1},\dots,y_r)$$
o\ \ $z_i = (x_i,y_i)$. Ainsi $v \in L_r(Z;F)$.

\

Soit $g(z) = v(z,\dots,z)$. Cela d\Ž finit un polyn\™ me $g
\in F[Z]_r$ et l'on a
$$\hat g(x,y) = g(x+y) = u(x,\dots,x,y,\dots,y) =
f(x,y).\qed$$
Cela ach\ ve la quatri\ me identification. La cinqui\ me
d\Ž coule des pr\Ž c\Ž dentes.

\

\

Bien entendu, ce que l'on vient de faire pour \bf deux \rm
variables se g\Ž n\Ž ralise au cas d'un nombre fini
quelconque de variables, celui des s\Ž ries formelles \bf
multiples\rm.

\

\su{1.6. Sur la notation $f(x)$} 

\

Soit $f \in F[E]$ un polyn\™
me homog\ ne d\Ž termin\Ž \ par $u \in L_m(E;F)$. Pour
chaque $x \in E$, l'\Ž l\Ž ment $f(x) = u(x,\dots,x)$ est un
vecteur de l'espace $F$. Plus g\Ž n\Ž ralement, lorsque $f
\in F[E]$ est un polyn\™ me quelconque, homog\ ne ou
non, pour chaque vecteur
$x$ de $E$, l'\Ž l\Ž ment $f(x)$ est encore un vecteur de
$F$. Autrement dit, $f$ d\Ž finit une application (polyn\™
miale)
$f : E \to F$.

\

Il n'en va plus n\Ž cessairement de m\ me pour une \bf
s\Ž rie \rm formelle 
$f \in F[[E]]$ quelconque. La notation $f(x)$ ne d\Ž signe
alors pas toujours un vecteur de $F$. Elle est pourtant
commode et on peut s'en servir (en abusant un peu) dans les
questions de composition ou de
\bf substitution\rm. En d\Ž signant encore cette s\Ž rie
par la notation $f : E
\to F$, on garde \sl m\Ž taphoriquement \rm l'id\Ž e d'une
application $f$ tout en sachant que $f(x)$ n'est pas un
vecteur de $F$. 

\

\su{1.7. Substitution}

\

Etant donn\Ž es des s\Ž ries formelles 
$$f : X \to E \ , \ g : E \to F,$$
on va montrer comment la compos\Ž e $g(f)$ peut avoir un
sens dans certains cas et qu'elle consiste \ˆ \ \bf
substituer
\rm
$f$ \ˆ \ $y$ dans $g(y)$. On fera cela par \Ž tapes.

\

\su{Substituter des polyn\™ mes homog\ nes dans une
application multilin\ aire} 

\

On commence par se donner 
$$u_i \in L_{n(i)}(X;E) \ , \ i = 1,\dots,m  \ , \ r = n(1)
+ \dots + n(m) \ , \ v \in
L_m(E;F).$$ 
L'application compos\Ž e
$w = v(u_1,\dots,u_m)$ appartient alors \ˆ \ $L_r(X;F)$. On
d\Ž signe par $f_1,\dots,f_m,h,$ les
polyn\™ mes homog\ nes d\Ž termin\Ž s, respectivement, par
$u_1,\dots,u_m,w,$ dont les degr\Ž s respectifs sont
$n(1),\dots,n(m),r$. Le
polyn\™ me $h$ n'est autre que le compos\Ž \
$v(f_1,\dots,f_m)$, on le voit imm\Ž diatement.

\

\su{Substituter des s\Ž ries formelles dans une
application multilin\ aire} 

\

Plus g\Ž n\Ž ralement, on se donne $m$ s\Ž ries formelles
$f_1,\dots,f_m,$ o\ \ $f_i \in E[[X]$. On d\Ž finit la s\Ž
rie formelle $v(f_1,\dots,f_m) = \sum_s h_s$ \ˆ \ l'aide de
ses composantes
$$h_s = \sum_{n(1) + \dots + n(m) = s}
v(f_{1,n(1)},\dots,f_{m,n(m)})$$
o\ \ $f_{i,j}$ d\Ž signe la composante de degr\Ž \ $j$ de
$f_i$.

\

\su{Substituer une s\Ž rie formelle dans un polyn\™ me}

\

Soit $f \in E[[X]]$ une s\Ž rie formelle et $g \in F[E]$ un
polyn\™ me. Lorsque $g$ est un polyn\™ me homog\ ne d\Ž
termin\Ž
\ par
$v
\in L_m(E;F)$, on d\Ž finit la s\Ž rie compos\Ž e par $g(f) =
v(f,\dots,f)$ qui ne d\Ž pend pas du choix de $v$, comme on
le v\Ž rifie. Lorsque $g(x) = \sum_{m \leq r}
v_m(x,\dots,x)$ est un polyn\™ me quelconque, avec $v_m
\in L_m(E;F)$, on
\Ž crit
$$g(f) = \sum_{m \leq r} v_m(f,\dots,f),$$
cette somme \Ž tant finie!

\

\su{Substituer une s\Ž rie formelle dans une autre}

\

Soient 
$$f : X \to E \ , \ g : E \to F$$
des s\Ž ries formelles. Pour chaque composante $g_s$ de $g$,
on consid\ re la s\Ž rie formelle $g_s(f)$. La famille
$(g_s(f))_s$ de ces s\Ž ries formelles  est sommable
uniquement dans les deux cas suivants, o\ \ l'on peut donc
\Ž crire $g(f) =
\sum_s g_s(f)$ :

\

\roster

\item lorsque $g$ est un polyn\™ me car il s'agit alors d'une
famille finie;

\

\item lorsque $g$ n'est pas un plolyn\™ me et que la
composante $f_0$ de degr\Ž \ $0$ de $f$ est nulle.
Sinon, les
$g_s(f_0)$ formeraient une famille infinie de constantes non
nulles.

\endroster

\

\su{1.8. Restriction et extension des coefficients et des
variables}

\

Chaque application
lin\Ž aire $p : Y \to X$ est un polyn\™ me homog\ ne
$p \in X[Y]_1$ de degr\Ž \ $1$ et induit, par
\sl fonctorialit\Ž
\rm,  une application lin\Ž aire
\sl naturelle \rm de $Y^m$ dans $X^m$ et, partant, des
applications lin\Ž aires 
$$L_m(X;E) \to L_m(Y;E) \ , \ E[X]_m \to E[Y]_m \ , \
\check p : E[[X]]
\to E[[Y]].$$
L'application $\check p$ respecte la \sl graduation \rm et,
pour
$f\in E[[X]]$, on a 
$\check p(f) = f(p)$  o\
\ l'on a substitu\Ž \ $p$ dans la s\Ž rie formelle $f$. 

\

De m\ me, \ˆ \ chaque application lin\Ž aire $s : E \to F$
correspondent des applications lin\Ž aires \sl naturelles \rm
$$L_m(X;E) \to L_m(X;F) \ , \ E[X]_m \to F[X]_m \ , \
\hat s : E[[X]]
\to F[[X]].$$
L'application $\hat s$ respecte aussi la \sl graduation \rm
et, pour
$f
\in E[[X]]$, on a
$\hat s(f) = s(f)$ o\
\ l'on a substitu\Ž \ la s\Ž rie formelle $f$ dans $s$.

\

En un certain sens, on peut ainsi dire que le foncteur
$E[[X]]$ est \bf contravariant
\rm en $X$ et \bf covariant \rm en $E$. On peut, \Ž
galement, voir les choses comme ceci : une extension $p : Y
\to X$ des variables induit une projection $\check p$ de
$E[[X]]$ dans
$E[[Y]]$; tandis qu'une extension $s : E \to F$ des
coefficients induit une extension $\hat s$ de $E[[X]]$ \ˆ \
$F[[X]]$.

\

Les deux applications $\check p$ et $\hat s$ commutent et
leurs effets conjugu\Ž s induisent une application lin\Ž
aire $q : E[[X]] \to \F[[Y]]$ o\ \ $q = \check p \circ \hat
s = \hat s \circ \check p$.

\

On pourrait examiner ce foncteur de plus pr\ s, afin
de voir comment se comportent les injections et les
surjections. Nous ne le ferons pas ici car nous n'en avons
pas l'usage. On se contentera de faire la seule
observation suivante.

\

Lorsque l'on a
$$E = X \ , \ F = Y \ , \ p \circ s = \roman{id}_E =
\roman{id}_X,$$
l'application $q$ est injective et plonge ainsi $X[[X]]$
dans $Y[[Y]]$.

\

\head 2. Sym\Ž trisation\endhead

\

\eightpoint Etant donn\Ž e une fonction de plusieurs
variables, $f(x_1,\dots,x_m)$, \ˆ \ valeurs dans un groupe
commutatif quelconque, une d\Ž marche habituelle de \bf
sym\Ž rtrisation\rm, classique,
 consiste
\ˆ
\ lui associer la fonction sym\Ž trique $g(x_1,\dots,x_m)$
obtenue en faisant la somme de toutes les valeurs
$f(x_{s(1)},\dots,x_{s(m)})$ o\ \  $s$ parcourt
l'ensemble des permuations de l'ensemble
fini $\{1,\dots,m\}$. On peut g\Ž n\Ž raliser quelque peu
cette d\Ž marche en d\Ž finissant la $(p_1,\dots,p_r)$-sym\Ž
tris\Ž e $\tilde f(z_1,\dots,z_r;p_1,\dots,p_r)$ de $f$
comme \Ž tant la somme de toutes les valeurs
$f(z_{s(1)},\dots,z_{s(m)})$ o\
\  $s$ parcourt l'ensemble des permuations \bf avec r\Ž p\Ž
titions
\rm de
$\{1,\dots,r\}$ o\ \ chaque $j$ est r\Ž p\Ž t\Ž \ $p_j$
fois. On explicite, ci-dessous, le cas particulier des
applications
$f$ multilin\Ž aires, pour les
besoins de la suite, notamment dans la d\Ž finition de la d\Ž
rivation suivant les s\Ž ries formelles, dans le paragraphe
\bf 3\rm.

\tenpoint

\

\su{2.1. Le $p$-sym\Ž tris\Ž \ d'une application $m$-lin\Ž
aire}

\

Pour chaque $u \in L_m(E;F) \ , \ z \in E^r \ , \ p =
(p_1,\dots,p_r) \in \Z^r$, on d\Ž signe par $\tilde u(z;p)$
la somme de tous les termes de la forme
$u(z_{s(1)},\dots,z_{s(m)})$ o\ \ $s$ est une permutation
avec r\Ž p\Ž titions de $\{1,\dots,r\}$ et dans laquelle $j$
est r\Ž p\Ž t\Ž \ exactement $p_j$ fois (pour $j =
1,\dots,r$). S'il n'existe aucun terme de cette forme, on
convient que $\tilde u(z;p) = 0$.

\

\su{Remarque} On observera qu'avec cette convention, on a
$\tilde u(z;p) = 0$ si $p \in \Z^r \setminus \N^r$ ou si $m
\neq p_1 + \dots + p_r$.

\

\su{2.2. Cas particulier} 

\

Pour d\Ž finir la d\Ž rivation suivant une s\Ž rie formelle,
nous nous servirons, ci-dessous, du cas particulier de la
sym\Ž trisation o\ \ $r = 2 \ , \ p_1 = 1 \ , \ p_2 = m-1$.

\

Dans ce cas, $\tilde u((t,x);(1,m-1))$ est la somme de tous
termes de la forme $u(x_1,\dots,x_m)$ o\ \ tous les $x_i$
sont \Ž gaux \ˆ \ $x$ sauf \bf un \rm qui est \Ž gal \ˆ \
$t$. Bien entendu, si $m = 0$ alors cette somme est nulle.

\

On va \Ž tablir le r\Ž sultat suivant.

\

\su{2.3. Th\Ž or\ me} \sl $<1>$. Si des \Ž l\Ž ments $u$ et
$v$  de $L_m(E;F)$ d\Ž terminent le m\ me polyn\™ me
homog\ ne, alors $\tilde u = \tilde v$.\rm

\

Autrement dit, si
$$u(x,\dots,x) = v(x,\dots,x) \ \ \text{pour tout} \ \ x \in
E,$$
alors
$$\tilde u(z;p) = \tilde v(z;p) \ \ \text{pour tous} \ \ z
\in E^r \ , \ p \in \Z^r.$$

\

Pour \Ž tablir ce th\Ž or\ me, on s'appuiera sur les deux
r\Ž sultats auxiliaires suivants.

\

\su{2.4. Lemme} Soient $b_0,b_1,\dots,b_m,$ une suite de
vecteurs dans $F$ et
$$g(t) = b_0 + t b_1 + \dots + t^m b_m \ \ \text{pour chaque}
\ \ t \in K.$$
On suppose que $g(t) = 0$ pour tout $t \in K$. Alors
$$b_0 = b_1 = \dots = b_m = 0.$$

\

\noi \bf En effet\rm, on consid\ re $(m+1)$ \Ž l\Ž ments
distinsts $t_0,t_1,\dots,t_m,$ dans $K$. On a
$$g(t_i) = 0 \ , \ i = 0,1,\dots,m, \ \text{par hypoth\
se}.$$
Autrement dit,
$$b_0 + t_0 b_1 + t_0^2 b_2 + \dots + t_0^m b_m = 0$$
$$b_0 + t_1 b_1 + t_1^2 b_2 + \dots + t_1^m b_m = 0$$
$$\vdots$$
$$b_0 + t_m b_1 + t_m^2 b_2 + \dots + t_m^m b_m = 0$$
et la matrice $M = (t^j_i)$ est (un Vandermonde) \bf
inversible\rm. D'o\ \ le r\Ž sultat.\qed

\

\su{Remarque} Bien entendu, le lemme est encore vrai lorsque
le corps $K$, sans \ tre infini, poss\ de au moins $(m+1)$
\Ž l\Ž ments distincts.

\

\su{2.5. Corollaire} $<2>$.

\

Soit $g(x_1,\dots,x_r)$ un polyn\™ me \ˆ \
coefficients dans $F$ de degr\Ž \ $m$ au plus :
$$g(x) = \sum_{p \in \N^r \ , \ |p| \leq m} x^p b_p.$$
On suppose que $g(x) = 0$ pour tout $x \in K^r$. Alors tous
les coefficients $b_p$ du polyn\™ me $g$ sont \bf nuls\rm.

\

\noi \bf En effet\rm, on raisonne par r\Ž currence sur
l'entier $r$ en utilisant le lemme pr\Ž c\Ž dent :

\noi $g(x_1,\dots,x_{r-1},t)$ est un polyn\™ me en
$(x_1,\dots,x_{r-1})$ de degr\Ž \ $\leq m$ \ˆ \ coefficients
dans $F[t]$.\qed

\

\su{2.6. Remarque} Soient $u \in L_m(E;F)$ et $z \in E^r$.
Voici une mani\ re commode de tenir serr\Ž s en une m\ me
formule tous les termes $\tilde u(z;p)$ pour $p \in \Z^r$.

\

Pour chaque $t = (t_1,t_2,\dots,t_r) \in K^r$, on consid\ re
$$x(t) = t_1 z_1 + t_2 z_2 + \dots + t_r z_r$$
et
$$f(t) = u(x(t),\dots,x(t)).$$
On a alors
$$f(t) = \sum_{p \in \N^r \ , \ |p| = m} t^p \tilde u(z;p)$$
en notation "multiindicielle" et o\ , par d\Ž finition, $|p|
= p_1 +
\dots + p_r$.

\

On se souvient, en effet, que $\tilde u(z;p) = 0$ dans
chacun des deux cas suivants :

\

(1) si $p \in \Z^r \setminus \N^r$,

\

(2) si$|p| \neq m$.

\

\su{2.7. D\Ž monstration du th\Ž or\ me 2.3} En consid\Ž
rant la diff\Ž rence $u-v$, il suffit d'\Ž tablir le r\Ž
sultat dans le cas o\ \ $v = 0$.

\

On suppose donc donn\Ž \ $u \in L_m(E;F)$ tel que
$u(x,\dots,x) = 0$ pour tout $x \in E$. On vient de montrer
que
$$\tilde u(z;p) = 0 \ \ \text{pour tous} \ \ z \in E^r \ , \
p
\in \N^r.$$
Or, en reprenant les notations de la remarque pr\Ž c\Ž
dente, on a
$$f(t) = u(x(t),\dots,x(t)) = 0 \ \ \text{pour tout} \ \ t
\in K.$$
On applique alors le \bf corollaire 2.5 \rm au polyn\™ me
homog\ ne $f(t)$ de degr\Ž \ $m$ en $t$, \ˆ \ coefficients
dans $F$.\qed 

\

\su{2.8. Remarque} Lorsque $m!$ est inversible dans $K$, \ˆ
\  chaque \Ž l\Ž ment $u
\in L_n(E;F)$ est associ\Ž \ son \bf sym\Ž
tris\Ž \ \rm 
$u^* \in L_m(E;F)$ que l'on d\Ž finit comme suit :
$$u^*(x_1,\dots,x_m) = \frac {1}{m!} \sum_s
u(x_{s(1)},\dots,x_{s(m)}),$$
la somme \Ž tant prise pour toutes les permutations $s$ de
l'ensemble $\{1,\dots,m\}$. Ce sym\Ž tris\Ž \ est une
application
$m$-lin\Ž aire \bf sym\Ž trique \rm de $E$ dans $F$. Lorsque
$f \in F[E]_m$ est d\Ž termin\Ž \ par $u$, il est clair que
$f$ est
\Ž galement d\Ž termin\Ž \ par ce sym\Ž tris\Ž \ $u^*$, l'\bf
unique \rm \Ž l\Ž ment sym\Ž trique de $L_m(E;F)$ qui d\Ž
termine $f$.

\

\head 3. D\Ž rivation suivant une s\Ž rie formelle\endhead

\

\eightpoint Au sens classique, la d\Ž riv\Ž e d'une s\Ž rie
formelle $f \in K[[x]]$, $f(x) = a_0 + a_1x + a_2x^2 + \dots
+ a_mx^m + \cdots$, est la s\Ž rie formelle $f'(x) = a_1 +
2a_2x +
\dots + ma_mx^{m-1} + \cdots\cdot$ Autrement dit, $f'(x)$
est le terme constant en $t$ de la s\Ž rie formelle $(f(x +
t) - f(x))/t \in K[[x,t]]$. On peut \Ž tendre, l\Ž g\
rement, cette d\Ž finition comme suit. Etant donn\Ž e une
autre s\Ž rie formelle
$\xi \in K[[x]]$, on appelle d\Ž riv\Ž e de $f$ suivant $\xi$
le terme constant en $t$ de la
s\Ž rie formelle $(f(x + t\xi) - f(x))/t$ : c'est la s\Ž rie
formelle $a_1\xi + 2\xi x + \dots m\xi x^{m-1} +
\cdots\cdot$  On d\Ž signe cette d\Ž riv\Ž e par
$\xi . f$. En particulier, la d\Ž riv\Ž e classique n'est
autre que $1.f = f'$, la d\Ž riv\Ž e de $f$ suivant la s\Ž
rie constante $1$. Ces consid\Ž rations se g\Ž n\Ž ralisent
au cas de s\Ž ries formelles $f \in F[[X]]$ et $\xi \in
X[[X]]$.

\tenpoint

\

\centerline{\bf Pour simplifier, on \Ž crira $S(X)$ au
lieu de $X[[X]]$, et $S_m(X)$ au lieu de $X[[X]]_m$\rm.}

\

\su{3.1. D\Ž riv\Ž e suivant un polyn\™ me homog\ ne}

\

Soit $\xi \in S_r(X)$ un polyn\™ me homog\ ne de degr\Ž \
$r$ \ˆ \ variables \bf et \rm coefficients dans $X$.

\

On va lui associer une \bf d\Ž rivation \rm qui, \ˆ \ chaque
s\Ž rie formelle $f \in F[[X]]$, fait correspondre une s\Ž
rie formelle $\xi.f \in F[[X]]$, la \bf d\Ž riv\Ž e de $f$
suivant $\xi$.\rm

\

Pour cela, on proc\ de par \Ž tapes.

\

\su{3.1.1} On commence par le cas o\ \ $f$ est un polyn\™
me homog\ ne. On suppose donc que $f \in F[X]_m$ est d\Ž
termin\Ž \ par $u \in L_m(X;F)$. Pour chaque $x \in E$, on
consid\ re
$$g(x) = \tilde u(\xi (x),x);(1,m-1))$$
voir ci-dessus, cas particulier, \bf 2.2\rm. 

\

Cela d\Ž finit une application $g : X \to F$. On commence
par observer que $g$ ne d\Ž pend pas du choix de $u$,
d'apr\ s le \bf th\Ž or\ me 2.3\rm. 

\

On d\Ž signera
l'application $g$ ainsi construite par $\xi.f$ et on
l'appellera \bf d\Ž riv\Ž e de $f$ suivant $\xi$\rm.

\

Autrement dit,
$$(\xi.f)(x) = u(\xi(x),x,\dots,x) + u(x,\xi(x),x,\dots,x) +
\dots + u(x,x,\dots,x,\xi(x)).$$

\

\su{ 3.1.2. Proposition} $<3>$.

\

Soient $\xi \in S_r(X)$ et $f \in F[X]_m$. Alors
$$\xi .f \in F[X]_s \ \ \text{o\ \ } \ \ s = r + m - 1,$$
avec la convention, naturellement, que $F[X]_s = \{0\}$ pour
$s < 0$.

\

\noi \bf En effet\rm, soient $u \in L_m(X;F)$ et $v \in
L_r(X;X)$ tels que
$$f(x) = u(x,\dots,x) \ \ \text{et} \ \ \xi (x) =
v(x,\dots,x).$$
Pour chaque $i=1,2,\dots,m,$ on d\Ž finit une application
$w_i : X^s \to F$ comme suit 
$$w_i(x_1,\dots,x_{i-1},x_{i+1},\dots,x_m,y_1,\dots,y_r) =
u(x_1,\dots,x_{i-1},v(y_1,\dots,y_r),x_{i+1},\dots,x_m).$$
On consid\ re ensuite $w = w_1 + w_2 + \dots + w_m$. Alors
$$(\xi .f)(x) = w(x,\dots,x) \ \ \text{et} \ \ w \in
L_s(X;F).\qed$$

\

\su{3.1.3. Le cas o\ \ $f$ est une s\Ž rie formelle} Lorsque
$f \in F[[X]]$, on d\Ž finit $\xi .f$ \ˆ \ l'aide des d\Ž
riv\Ž es des composantes homog\ nes de $f$. Ainsi, pour
$$f = \sum_m f_m \ \ \text{o\ \ } \ \ f_m \in F[X]_m,$$
on pose
$$\xi .f = \sum_m \xi .f_m.$$

\

\su{3.2. D\Ž riv\Ž e suivant une s\Ž rie formelle}

\

Soient 
$$\xi \in S(X) \ \ \text{et} \ \ f \in F[[X]],$$
$$\text{o\ \ } \ \ \xi = \sum \xi_r \ \ \text{et} \ \ f =
\sum f_m,$$
les $\xi_r$ et les $f_m$ \Ž tant les composantes homog\ nes.

\

G\Ž n\Ž ralisant la d\Ž finition ci-dessus, on d\Ž signera
par $\xi .f$ la s\Ž rie formelle $\sum g_s$ o\ , pour $s
\geq 0$, on a pos\Ž
$$g_s = \sum_{r+m-1=S} \xi_r .f_m.$$

\

L'application de $S(X) \times F[[X]]$ dans $F[[X]]$ qui au
couple $(\xi,f)$ fait correspondre la d\Ž riv\Ž e $\xi.f$
est une application \bf bilin\Ž aire\rm.

\

\su{3.3. Remarque} Lorsque le corps $K$ est de \bf caract\Ž
ristique nulle\rm, toutes les factorielles $m!$ y sont
inversibles et chaque polyn\™ me homog\ ne $f_m \in T[X]_m$
est alors d\Ž termin\Ž \ par une application $m$-lin\Ž
aires \bf sym\Ž trique \rm $u_m = u^*_m \in L_m(X;F)$ (voir
la
\bf Remarque 2.8\rm). Toute s\Ž rie formelle $f \in
F[[X]]$ a un d\Ž veloppement canonique
$$f(x) = \sum_m u_m(x,x,\dots,x),$$
o\ \ les $u_m$ sont sym\Ž triques. Pour la d\Ž riv\Ž e de
$f$ suivant $\xi$, on a alors la formule particu\-li\ rement
simple et suggestive suivante :
$$(\xi . f)(x) = \sum_m m \ u_m(\xi(x),x,\dots,x).$$

\

\su{3.4. Plongement et d\Ž rivation}

\

On reprend les notations du \bf 1.5 \rm pour les s\Ž ries
doubles :
$$T = F[[X]] \ , \ T_m = F[[X]]_m \ , \ Z = X \times Y.$$
On a identifi\Ž \ $T[[Y]]$ \ˆ \ $F[[Z]]$. De m\ me,
$S(X)[[Y]]$ est identifi\Ž \ \ˆ \ $X[[Z]]$ lequel
est plong\Ž \ dans $S(Z)$.

\

Soient
$$f \in T_m[Y]_n \subset F[[Z]] \ , \ \xi \in S_p(X)[Y]_q
\subset S[[Z]].$$
Pour $y \in Y$, on a
$$f(y) \in T_m = F[X]_m \ , \ \xi(y) \in S_p(X) \subset
S(X).$$
Soit
$$g(y) = \xi(y).f(y)$$
la d\Ž riv\Ž e du polyn\™ me $f(y)$ suivant $\xi(y)$; de
sorte que l'on a $g(y) \in F[X]{m+p-1}$. Soit alors $h =
\xi.f$ la d\Ž riv\Ž e de $f \in F[[Z]]$ suivant $\xi \in
S(Z)$; de sorte que l'on a $h \in F[[Z]]$.

\

\su{3.5. On montre alors que l'on a $h(x,y) = g(y)(x)$}

\

Autrement dit, \bf calculer la d\Ž riv\Ž e $\xi.f$ \lq\lq
point par point" pour chaque $y$, revient \ˆ \ calculer
cette d\Ž riv\Ž e \lq\lq globalement" pour \sl la variable
dans $Z$\rm.

\

\su{D\Ž monstration} Si $u \in L_n(Y;L_m(X;F)$ d\Ž termine
$f$, on consid\ re $v \in L_r(Z;F)$, avec $r = m + n$, d\Ž
finie par
$$v(z_1,\dots,z_r) = u(y_{m+1},\dots,y_r)(x_1,\dots,x_m) \
, \ \text{o\ } \ \ z_i = (x_i,y_i),$$
de sorte que $v$ d\Ž termine $f$ comme \Ž l\Ž ment de
$F[[Z]]$ car
$$f(z) = f(y)(x) = u(y,\dots,y)(x,\dots,x) = v(z,\dots,z) \
, \ \text{o\ } \ \ z = (x,y).$$

\noi \bf{On calcule}
$$g(y)(x) = (\xi(y).f(y))(x) = \sum_i
u(y,\dots,y)(x,\dots,\xi(x)\urcorner^i,\dots,x).$$
Ensuite $\xi$ comme \Ž l\Ž ment de $S(z)$ :
$$\xi(z) = \xi(x,y) = (\xi(y)(x),0).$$
Enfin,
$$h(z) = h(x,y) = (\xi.f)(z) = \sum
v(z,\dots,\xi(z),\dots,z) =$$
$$\sum_i u(y,\dots,y)(x,\dots,\xi(x)\urcorner^i,\dots,x) =
g(y)(x).\qed$$\rm

\

\su{3.6. Deux d\Ž finitions}

\

On se donne une partie $N \subset X$ de l'espace vectoriel
$X$ ainsi que deux s\Ž ries formelles $\xi \in S(X)$ et $f
\in F[[X]]$.

\

On dira que $\xi$ \bf prend ses valeurs \rm dans $N$
lorsque, pour chacune de ses composantes homog\ nes $\xi_m$
et pour chaque $x \in X$, on a
$$\xi_m(x) \in N.$$
On dira que \bf $f$ ne d\Ž pend pas de $N$ \rm lorsque
chacune de ses composantes homog\ nes $f_n$ peut \ tre d\Ž
termin\Ž e par une application $n$-lin\Ž aire $u \in
L_n(X;F)$ ayant la propri\Ž t\Ž \ suivante :
$$u(x_1,\dots,x_n) = 0 \ \ \text{d\ s que l'un des } \ 
x_i \ \ \text{appartient \ˆ } \ \ N.$$
En particulier, on observera donc que, lorsque $f$ ne d\Ž
pend pas de
$N$, on a 
$$f(x) = 0 \ \ \text{pour tout} \ \ x \in N.$$

\

\su{3.7. Lemme} \sl Si $\xi$ prend ses valeurs dans $N$ et
$f$ ne d\Ž pend pas de $N$, la d\Ž riv\Ž e $\xi.f$ est
nulle\rm.

\

\noi \bf En effet\rm, il suffit de faire la d\Ž monstration
pour les composantes homog\ nes $\xi_m$ et $f_n$. Or,
$$(\xi_m.f_n)(x) = \sum
u(x,\dots,\xi_m(x)\urcorner^i,\dots,x),$$
et chacun des termes de cette somme est nul.\qed

\

\su{3.8. D\Ž rivation d'une s\Ž rie compos\Ž e} D'un seul
mot, on signalera encore cette simple formule dont on se
servira plus loin. Lorsque $f$ est lin\Ž aire, autrement dit,
lorsque $f$ est un polyn\™ me homog\ ne de degr\Ž \ $1$, on
a
$$\xi.(f(g)) = f(\xi . g).$$

\head 4. L'alg\ bre de Lie $S(X)$\endhead

\

\eightpoint L'ensemble de s\Ž ries formelles $S(X)$ muni de
l'op\Ž ration de d\Ž rivation $(\xi,\eta) \mapsto \xi .
\eta$ est une alg\ bre, mais elle n'est \bf pas associative
\rm (sauf cas trivial). Pourtant, fait inattendu, son crochet
$[\xi,\eta] = \xi . \eta - \eta . \xi$ v\Ž rifie
l'identit\Ž \ de Jacobi et en fait une alg\ bre de Lie.

\tenpoint

\

Etant donn\Ž es des s\Ž ries formelle $\xi$ et $\eta$, \Ž l\Ž
ments de $S(X)$, on peut consid\Ž rer $\xi.\eta$, la d\Ž
riv\Ž ede $\eta$ suivant $\xi$, puis $\eta.\xi$, la d\Ž
riv\Ž e de $\xi$ suivant $\eta$. Ce sont deux s\Ž ries
formelles \Ž l\Ž ments de $S(X)$.

\

On pose alors
$$[\xi,\eta] = \xi.\eta - \eta.\xi.$$
On appelle cette op\Ž ration le \bf crochet \rm sur $S(X)$.

\

On obtient alors le r\Ž sultat (inattendu !) suivant.

\

\su{4.1. Th\Ž or\ me} $<4>$. \sl L'espace vectoriel $S(X)$
muni du crochet est une alg\ bre de Lie\rm.

\

On proc\Ž dera par \Ž tapes pour \Ž tablir ce th\Ž or\ me.

\

\su{4.2. Lemme} Soient 
$$\xi \in S(X) \ , \ \eta \in S(X) \ ,
\ f \in F[[X]].$$
Alors
$$\xi.(\eta.f) - \eta.(\xi.f) = (\xi.\eta - \eta.\xi).f.$$
Autrement dit,
$$[\xi,\eta].f = \xi.(\eta.f) - \eta.(\xi.f).$$

\

\noi \bf En effet\rm, par bilin\Ž arit\Ž , on se ram\ ne au
cas des polyn\™ mes homog\ nes.

\

Soient alors 
$$\xi \in S_r(X) \ , \ \eta \in S_s(X) \ , \ f \in F[X]_m.$$
On suppose que $f,\xi,\eta,$ sont d\Ž termin\Ž s
respectivement par $u,v,w$.

\

Calculons, avec des notations qui s'expliquent d'elles-m\
mes :
$$(\xi.f)(x) = \sum_i u(x,\dots,\xi(x)\urcorner^i,\dots,x)$$
$$\eta.(\xi.f)(x) = \sum_{i \neq j}
u(x,\dots,\xi(x)\urcorner^i,\dots,\eta(x)\urcorner^j,
\dots,x) \ +$$
$$+ \sum_{i , j}
u(x,\dots,v(x,\dots,\eta(x)\urcorner^j,
\dots,x)\urcorner^i,\dots,x).$$
De m\ me,
$$\xi.(\eta.f)(x) = \sum_{i \neq j}
u(x,\dots,\xi(x)\urcorner^i,\dots,\eta(x)\urcorner^j,
\dots,x)$$
$$+ \sum_{i , j}
u(x,\dots,w(x,\dots,\xi(x)\urcorner^j,
\dots,x)\urcorner^i,\dots,x).$$
De sorte que
$$(\xi.(\eta.f) - \eta.(\xi.f))(x) = \sum_i
u(x,\dots,((\xi.\eta)(x) -
(\eta.\xi)(x))\urcorner^i,\dots,x) = ((\xi.\eta -
\eta.\xi).f)(x).\qed$$

\

\su{4.3. Remarque} Soient $\xi,\eta,\zeta,$ des s\Ž ries
formelles, \Ž l\Ž ments de $S(X)$. Les deux s\Ž ries formelle
$\xi.(\eta.\zeta)$ et $(\xi.\eta).\zeta$ ne sont pas n\Ž
cessairement \Ž gales.

\

Donnons un contre-exemple tr\ s simple! Pour 
$$X = K \ , \ \xi(x) = \zeta(x) = x^2 \ , \ \eta(x) = x,$$
on a
$$(\xi.\eta)(x) = x^2 = \xi(x),$$
$$(\xi.\eta).\zeta = \xi.\xi \ , \ (\xi.\xi)(x) = 2x^3,$$
$$(\eta.\zeta)(x) = 2x^2,$$
$$\xi.(\eta.\zeta)(x) = 4x^3.$$
\bf Autrement dit, l'op\Ž ration de d\Ž rivation n'est pas
associative\rm.

\

\su{4.4. Lemme [\rm identit\Ž \ de Jacobi\rm]} Le crochet de
$S(X)$ v\Ž rifie l'identit\Ž \ de Jacobi :
$$[a,[b,c]] = [[a,b],c] -[[a,c],b].$$

\

\noi \bf En effet\rm,  en utilisant le lemme pr\Ž c\Ž dent
\bf 4.2\rm, on a
$$[[a,b],c] - [[a,c],b] =$$
$$(ab - ba)c - c(ab - ba) - (ac - ca)b + b(ac - ca) =$$
$$a(bc) - b(ac) - c(ab) + c(ba) - a(cb) + c(ab) + b(ac) -
b(ca) =$$
$$a(bc) + c(ba) - a(cb) - b(ca) =$$
$$a(bc - cb) + c(ba) - b(ca) =$$
$$a[b,c] - (bc - cb)a = a[b,c] - [b,c]a = [a,[b,c]].\qed$$

\

Cela d\Ž montre le \bf th\Ž or\ me 4.1\rm.

\

\head 5. Action d'une alg\ bre de Lie sur un espace
vectoriel
\endhead

\

\eightpoint \`A chaque espace vectoriel
$X$ est attach\Ž , naturellement, une alg\ bre de Lie
$S(X)$, comme on vient de le voir. Il est tout aussi naturel
de fair agir une alg\ bre de Lie donn\Ž e $A$ sur l'espace
vectoriel $X$ par le biais des homomorphismes $\rD : A \to
S(X)$ d'alg\ bres de Lie.

\tenpoint

\

On appellera \bf action  formelle \rm (\ˆ \ doite) de
l'alg\ bre de Lie $A$ sur l'espace vectoriel $X$ tout
homomorphisme d'alg\ bres de Lie $\rD : A \to S(X)$.

\

Ainsi, pour chaque $a \in A$, l'\Ž l\Ž ment correspondant
$\rD_a \in S(X)$ est une s\Ž rie formelle \ˆ \ variables et
coefficients dans $X$. Pour chaque $\lambda \in K \ , \ a \in
A \ , \ b \in A$, on aura donc
$$\rD_{\lambda a} = \lambda \rD_a \ , \ \rD_{a+b} = \rD_a +
\rD_b \ , \ \rD_{[a,b]} = [\rD_a,\rD_b].$$
De plus, pour toute s\Ž rie formelle $f \in F[[X]]$, la d\Ž
riv\Ž e $\rD_a f$ (alias $\rD_a.f$) de $f$ suivant $\rD_a$
est, elle-m\ me, une s\Ž rie formelle, \Ž l\Ž ment de
$F[[X]]$ \Ž galement. On a, aussi, la formule suivante qui
d\Ž coule du
\bf lemme 4.2 \rm :

$$\rD_{[a,b]} f = [\rD_a,\rD_b] f = (\rD_a \rD_b -
\rD_b \rD_a) f = \rD_a(\rD_b f) -
\rD_b(\rD_a f).$$

\

\head Prolongement canonique d'une action formelle\endhead

\

Soit $\rD : A \to S(X)$ une action formelle de l'alg\ bre
de Lie $A$ sur l'espace vectoriel $X$. Pour chaque espace
vectoriel $Y$, on d\Ž finira une struture d'alg\ bre de Lie
naturelle sur $A[[Y]]$ puis un prolongement canonique de
l'action $\rD$ en une action de l'alg\ bre de Lie $A[[Y]]$
sur l'espace vectoriel produit $X \times Y$.

\

\su{5.1. Crochet sur $A[[Y]]$}

\

Ce sera simplement le \lq\lq crochet ponctuel", ou encore
\lq\lq point par point", \bf h\Ž rit\Ž \ \rm de $A$.

\

Autrement dit, si $f$ et $g$
appartiennent \ˆ \ $A[[Y]]$ et $y \in Y$, on dira, de mani\
re imag\Ž e, que
$$[f,g](y) = [f(y),g(y)].$$
Plus pr\Ž cis\Ž ment, et techniquement : si $f \in A[Y]_n$
et $g \in A[Y]_r$, pour chaque $y \in Y$, on pose
$$[f,g](y) = [f(y),g(y)].$$ 
Cela d\Ž finit une application $[f,g] : Y \to A$.

\

\su{5.2. Lemme} $<6>$. Soient
$$f \in A[Y]_n \ \ \text{et} \ \ g \in A[Y]_r;$$
Alors
$$[f,g] \in A[Y]_{n+r}.$$
\bf Plus pr\Ž cis\Ž ment\rm, si $f$ et $g$ sont d\Ž termin\Ž
s respectivement par $u$ et $v$, on consid\ re l'application
$$w : Y^{n+r} \to A$$
d\Ž fine par
$$w(x_1,\dots,x_n,y_1,\dots,y_n) =
[u(x_1,\dots,x_n),v(y_1,\dots,y_n)] \ \ \text{pour le crochet
de} \ \ A.$$
Alors $w \in L_{n+r}(Y;A)$. On \Ž crira simplement $w =
[u,v]$. Enfin, $[f,g]$ est d\Ž termin\Ž \ par $[u,v]$\qed 

\

\su{5.3. L'alg\ bre de Lie $A[[Y]]$}

\

A pr\Ž sent, lorsque $f$ et $g$ sont des s\Ž ries formelles
\Ž l\Ž ments de $A[[Y]]$, o\ \
$$f = \sum f_n \ , \ g = \sum g_r,$$
on pose
$$[f,g]_s = \sum_{n+r=s} [f_n,g_r]$$
et on d\Ž finit $[f,g]$ comme \Ž tant la s\Ž rie formelle
$\sum [f,g]_s$.

\

Ce crochet sur $A[[Y]]$, \bf h\Ž rit\Ž \ \rm du crochet sur
$A$, est une application bilin\Ž aire antisym\Ž trique d\Ž
finie sur $A[[Y]]$.

\

\su{5.4. Proposition} $<7>$.
\sl Le crochet h\Ž rit\Ž \ de $A$ munit
l'espace vectoriel $A[[Y]]$ d'une structure d'alg\ \-bre de
Lie\rm.

\

\noi \bf En effet\rm, \lq\lq point par point", le crochet
v\Ž rifie l'identit\Ž \ de Jacobi.\qed

\

Comme annonc\Ž \ ci-dessus, on a donc muni l'espace $A[[Y]]$
d'une structure d'alg\ bre de Lie h\Ž rit\Ž e de celle de
$A$.

\

\su{5.5. Prolongement canonique d'une action formelle}

\

Soit $\rD : A \to S(X)$ une action formelle de l'alg\ bre de
Lie
$A$ sur l'espace vectoriel $X$. Soit $Y$ un espace vectoriel
quelconque. Consid\Ž rons l'espace vectoriel produit $Z = X
\times Y$ et l'alg\ bre de Lie $A[[Y]]$.

\

Soit $a \in A[[Y]]$. On peut \bf substituer \rm la s\Ž rie
formelle $a$ dans l'application lin\Ž aire $\rD$ (voir
ci-dessus au \bf 1.7\rm). D\Ž signons le r\Ž sultat de cette
substitution par $\rD_a$. Ainsi $\rD_a$ est une s\Ž rie
formelle \ˆ \ variables dans $Y$ et coefficients dans
$S(X)$. Autrement dit,
$$\rD_a \in S(X)[[Y]] = (X[[X]])[[Y]].$$
Ce dernier espace de s\Ž ries formelles est canoniquement
plong\Ž \ dans $X[[X \times Y]]$, lui-m\ me plong\Ž \ dans
$Z[[Z]] = S(Z)$.

\

Finalement, \ˆ \ chaque $a \in A[[Y]]$ correspond une s\Ž
rie formelle $\rD_a \in S(Z)$. Nous d\Ž signerons cette
application par la m\ me lettre
$$\rD : A[[Y]] \to S(Z).$$
C'est visiblement une application lin\Ž aire et qui prolonge
$\rD : A \to S(X)$.

\

\su{5.6. Proposition} 

\

\sl Soit $\rD : A \to S(X)$ une action
formelle de l'alg\ bre de Lie $A$ sur l'espace vectoriel
$X$.
Pour chaque espace vectoriel $Y$, le prolongement $\rD :
A[[Y]] \to S(X \times Y)$ est une action formelle de l'alg\
bre de Lie $A[[Y]]$ sur l'espace vectoriel produit $Z = X
\times Y$\rm.

\

\noi \bf En effet\rm, il suffit de prouver que
$$\rD_{[a,b]} = [\rD_a,\rD_b]$$
pour tout couple de s\Ž ries formelles $a$ et $b$ \Ž l\Ž
ments de $A[[Y]]$. Il suffit donc d'\Ž tablir cette \Ž
galit\Ž \ lorsque $a$ et $b$ sont des polyn\™ mes homog\
nes.

\

Soient 
$$a \in A[Y]_n \ \ \text{et} \ \ b \in A[Y]_r.$$
Alors
$$[a,b] \in A[Y]_{n+r} \ \ \text{et} \ \ [a,b](y) =
[a(y),b(y)] \ \ \text{pour chaque} \ \ y \in Y.$$
Comme $\rD : A \to S(X)$ est une action formelle, on a
$$\rD_{[a,b](y)} = \rD_{[a(y),b(y)]} =
[\rD_{a(y)},\rD_{b(y)}] \ \ \text{pour chaque} \ \ y \in Y.$$
D'un autre c\™ t\Ž , consid\Ž rons $[\rD_a,\rD_b] =
\rD_a.\rD_b - \rD_b.\rD_a$. C'est un \Ž l\Ž ment de $S(Z)$
que l'on peut calculer \lq\lq ponctuellment" comme on l'a
vu plus haut (au \bf 3.5\rm). Ce n'est autre que la s\Ž rie
formelle dans $S(Z)$ qui correspond au polyn\™ me homog\ ne
$h \in S(X)[Y]$ o\ \
$$h(y) = \rD_{a(y)}.\rD_{b(y)} - \rD_{b(y)}.\rD_{a(y)} =
[\rD_{a(y)},\rD_{b(y)}],$$
d'o\ \ l'\Ž galit\Ž \.\qed

\

\su{5.7. Sur le plongement de $S(X)[[Y]]$ dans S(Z)}

\

On reprend les notations du \bf 1.5 \rm pour les s\Ž ries
doubles; dans le cas o\ \ $F = X$. Ainsi :
$$T = S(X) \ , \ Z = X \times Y.$$
On a vu que $T[[Y]]$ s'identifie canoniquement \ˆ \ $X[[X
\times Y]] = X[[Z]]$ lequel se plonge dans $Z[[Z]] = S(Z)$.
D'o\ \ un plongement canonique
$$j : T[[Y]] \to S(Z).$$
On a muni $T = S(X)$ d'une structure d'alg\
bre de Lie (voir au \bf 3\rm, ci-dessus) dont $T[[Y]]$ h\Ž
rite
\lq\lq ponctuellement" (voir au \bf 5.4\rm, ci-dessus).
D'autre part $S(Z)$ est \Ž galement muni, intrins\Ž quement,
d'une stucture d'alg\Ž bre de Lie.

\

On \Ž tablit le r\Ž sultat \sl naturel \rm suivant :

\

\su{5.7.1. Th\Ž or\ me} \sl Le plongement $j : T[[Y]] \to
S(Z)$ est un plongement d'alg\ bres de Lie\rm. $<11>$.

\

Autrement dit, pour toutes s\Ž ries formelles $f$ et $g$ dans
$T[[Y]]$, on a
$$j([f,g]) = [j(f),j(g)].$$
\bf En effet\rm, on proc\ de en deux \Ž tapes. On commence
par supposer que
$$f \in T_m[X]_n \ , \ g \in T_p[X]_q.$$
Soit $h = [f,g]$, crochet dans $T[[Y]]$. Ainsi, pour chaque
$y \in Y$, on a
$$h(y) = [f(y),g(y)] = f(y).g(y) - g(y).f(y).$$
Ici, $f(y).g(y)$ est la d\Ž riv\Ž e (sur la variable dans
$X$) de $g(y)$ suivant $f(y)$. D'apr\ s le r\Ž sultat ant\Ž
rieur du \bf 3.5 \rm ci-dessus, on a
$$(f(y).g(y))(x) = (j(f).(j(g))(x,y),$$
o\ \ $j(f).j(g)$ est la d\Ž riv\Ž e de $j(g)$ suivant $j(f)$
pour
\sl la variable dans $Z$\rm. Donc
$$h(y)(x) = [j(f),j(g)](x,y),$$
avec le crochet $[f,g]$ de $S(Z)$, et
$$j(h) = [j(f),j(g)].$$
Dans le cas g\Ž n\Ž ral, on applique ce qui pr\Ž c\ de aux
composantes homog\ nes de $f$ et de $g$ et on se sert de la
double additivit\Ž \ du crochet.\qed

\

\su{5.7.2 Plongement de l'alg\ bre de Lie $S(X)$ dans
l'alg\ bre de Lie $S(Z)$} $<10>$.

\

Au passage, on peut observer que l'on obtient ainsi, en
particulier, un plogement \sl naturel \rm de l'alg\ bre de
Lie $S(X)$ dans l'alg\ bre de Lie $S(Z)$, sans nouvelle
argumentation. \bf En effet\rm, l'ensemble $T[Y]_0$ des
\sl constantes \rm de $T[[Y]]$ est une sous-alg\ bre de Lie
de $T[[Y]]$ isomorphe \ˆ \ $S(X)$.\qed

\

De m\ me, bien entendu, $S(Y)$ est \Ž galement plong\Ž e
dans
$S(Z)$.

\

\head 6. L'exemple originel \endhead

\

\eightpoint

Nous avons introduit les notions \sl formelles \rm pr\Ž c\Ž
dentes en nous inspirant de l'exemple de la notion de \sl
loi d'op\Ž ration infinit\Ž simale \rm selon Bourbaki.  Cet
exemple originel des lois d'op\Ž rations infinit\Ž
simales, ce prototype
\bf analytique\rm, nous a servi de mod\ le pour
l'introduction des notions purement \bf alg\Ž briques \rm de
d\Ž rivation et d'action formelles pour les alg\ bres de
Lie.

\tenpoint

\

On suppose, ici, que $K = \R$ ou $\C$ et on se donne
des espaces de Banach $X$
et $F$. 

\

Dans ce cas, on d\Ž signe par $\hat P(X;F)$ l'ensemble des
\bf s\Ž ries formelles \ˆ \ composantes continues \rm sur
$X$ \ˆ
\ valeurs dans $F$ (voir Bourbaki [1], p.88-89). Ainsi $\hat
P(X;F)$ est un sous-espace vectoriel de $F[[X]]$, qui lui
est \Ž gal lorsque la dimension de $X$ est \bf finie\rm.

\

\su{6.1. D\Ž rivation suivant un champ de vecteurs}

\

Soient $U$ un voisinage ouvert de $0$ dans
$X$ et
$f : U
\to F$ une application analytique de $U$ dans $F$, tous deux,
consid\Ž r\Ž s comme des vari\Ž t\Ž s analytiques. Soit
$\xi$ un champ analytique de vecteurs sur $U$. On sait
d\Ž finir l'application analytique $\xi (f) : U \to F$
(voir  Bourbaki [2], 8.2.2 et 8.2.3,
p.10) : c'est la fonction $x \mapsto <\xi,\rd_x f>(x)$
o\
\
$\rd_x f$ d\Ž signe la diff\Ž rentielle de $f$ au point $x$.

\

En particulier, l'injection canonique $h : U \to X$ est
analytique, de m\ me que l'appli\-cation $\xi (h) : U \to
X$.

\

Identifions chacune des applications $\xi (h)$,
$f$ et $\xi(f)$ aux s\Ž ries formelles qui les repr\Ž
sen\-tent  au voisinage de $0$,
$$\xi(h) \in \hat P(X;X) \ , \ f \in \hat P(X;F) \ , \ 
\xi(f) \in \hat P(X;F).$$ 

\

On v\Ž rifie que l'on a
$$\xi(f) = \xi(h).f.$$ 
Autrement dit, la s\Ž rie formelle $\xi(f)$ n'est autre
que la d\Ž riv\Ž e de la s\Ž rie formelle $f$ suivant la s\Ž
rie formelle
$\xi(h)$, telle que cette d\Ž riv\Ž e a \Ž t\Ž \ d\Ž finie
ci-dessus (au
\bf 3\rm).

\

\su{V\Ž rification} Il suffit de la faire dans le cas o\ \
$f$ est un polyn\™ me homog\ ne continu de degr\Ž \ $m$ sur
$X$ \ˆ \ valeurs dans $F$, autrement dit,
lorsqu'il existe une application multilin\Ž aires continue
$X^m \to F$ telle que $f(x) = u(x,\dots,x)$, autrement dit,
$f = u(h,\dots,h)$. Dans ce cas, pour calculer $\xi(f)$, on
se sert de la formule de d\Ž rivation des fonctions
multilin\Ž aires compos\Ž es (voir Bourbaki [2], 8.2.3, page
11) :
$$\xi(f) = \sum_i
u(x,\dots,\xi(h)\urcorner^i,\dots,x) = \xi(h).f.\qed$$ 

\

\su{6.2. Lois d'op\Ž rations infinit\Ž simales} 

\

On reprend le
voisinage ouvert
$U$ de $0$ dans $X$ et on se donne
une alg\ bre de Lie normable compl\ te $A$. Une \bf
loi d'op\Ž ration infinit\Ž simale \ˆ \ droite\rm,
analytique, de l'alg\ bre de Lie $A$ dans la vari\Ž t\Ž \
analytique $U$ est une application $a \mapsto \xi_a$, o\ \
$\xi_a$ est un champ de vecteurs sur $U$,
ayant les deux propri\Ž t\Ž s suivantes (voir
Bourbaki [4], p.139).

\

\noi \bf (i) \rm L'application $(a,x) \mapsto \xi_a(x)$ est
un morphisme analytique du fibr\Ž \ vectoriel trivial $A
\times U$ dans le fibr\Ž \ tangent $T(U)$ lequel s'identifie
au fibr\Ž \ trivial $U \times X$.

\

\noi \bf (ii) \rm On a $[\xi_a,\xi_b] = \xi_{[a,b]}$ quels
que soient $a$ et $b$ dans $A$.

\

En particulier, pour chaque $a \in A$, le champ de
vecteurs
$\xi_a$ est analytique. Bien entendu, l'injection canonique
$h : U \to X$ est analytique. On consid\ re l'application
$\xi_a(h) : U \to X$ (voir Bourbaki [2], 8.2.2 et 8.2.3,
p.10) : c'est la fonction $x \mapsto \rd_x h(\xi_a(x))$ o\ \
$\rd_x h$ d\Ž signe la diff\Ž rentielle de $h$ au point $x$.

\

Cette application $\xi_a(h)$ est analytique donc repr\Ž
sentable au voisinage de l'origine par une s\Ž rie formelle
(convergente) \ˆ \ composantes continues, c'est-\ˆ -dire par
un \Ž l\Ž ment de
$\hat P(X;X) \subset S(X)$, que nous d\Ž signerons par
$\rD_a$.

\

D'apr\Ž s ce qui a \Ž t\Ž \ dit, ci-dessus (au
\bf 6.1\rm), la d\Ž rivation suivant cette s\Ž rie
formelle $\rD_a$ op\ re comme suit :

\

Etant donn\Ž e une application $f : U \to F$, analytique au
voisinage de $0$, on a
$$\xi_a(f) = D_a.f.$$

On peut alors v\Ž rifier, simplement, que
\bf l'application
$\rD : A
\to S(X)$ ainsi d\Ž finie est, \rm{au sens que
nous lui avons donn\Ž }, \bf une action formelle de
l'alg\ bre de Lie $A$ sur l'espace vectoriel $X$\rm. 

\

On
dira que c'est l'action formelle
\bf d\Ž duite \rm de la loi d'op\Ž ration infinit\Ž simale
$\xi$ donn\Ž e.

\

\su{6.3. Th\Ž or\ me} $<5>$. \sl L'application $\rD : A \to
S(X)$ d\Ž duite d'une loi d'op\Ž ration infinit\Ž simale
$\xi$ donn\Ž e est une action formelle de $A$ sur $X$\rm.

\

\noi \bf En effet\rm, la seule chose qui pourrait ne pas \
tre tout \ˆ \ fait \sl claire\rm, c'est que le crochet
$[\rD_a,\rD_b]$ \bf calcul\Ž
\ dans l'alg\ bre de Lie
$S(X)$ \rm est \Ž gal \ˆ \ $\rD_{[a,b]}$.

\

On a 
$$\rD_{[a,b]} = \xi_{[a,b]}(h) = [\xi_a,\xi_b](h) =
\xi_a(\xi_b(h)) - \xi_b(\xi_a(h)) = \xi_a(\rD_b) -
\xi_b(\rD_a).$$ 
Or, $\xi_a(\rD_b) = \rD_a.\rD_b$, la d\Ž riv\Ž de la s\Ž rie
formelle $\rD_b$ suivant la s\Ž rie formelle
$\rD_a$. De sorte que
$$\rD_{[a,b]}= \rD_a.\rD_b - \rD_b.\rD_a =
[\rD_a,\rD_b].\qed$$

\

\head 7. Produits d'entrelacements\endhead

\

\eightpoint

On se donne une action formelle $\rd : B \to S(Y)$ de
l'alg\ bre de Lie $B$ sur $Y$. On introduit l'espace
vectoriel produit $W = A[[Y]] \times B$ et on  d\Ž
finit une structure d'alg\ bre de Lie sur $W$ que l'on
appellera \bf produit d'entrelacement de $B$ par $A$
relativement \ˆ \ l'action formelle $\rd$\rm. On d\Ž
signera cette alg\ bre de Lie par $W(A,B;\rd)$. Elle
se pr\Ž sentera comme un \bf produit semi-direct
\rm d'alg\ bres de Lie. A chaque action formelle $\rd : B
\to S(Y)$, un produit d'entrelacment $W(A,B;\rd)$ :
c'est la raison du pluriel.

\tenpoint

\

\su{7.1. L'alg\ bre de Lie $\Der(A : Y)$}

\

On va d\Ž signer par $\Der(A : Y)$ l'alg\ bre de Lie des d\Ž
rivations de l'alg\ bre de Lie $A[[Y]]$. Autrement dit,
$\Der(A : Y) = \frak d(A[[Y]])$.

\

Pour tout \Ž l\Ž ment $f \in S(Y)$ et tout $a \in A[[Y]]$, on
peut consid\Ž rer la d\Ž riv\Ž e $f.a$ de la s\Ž rie $a$
suivant la s\Ž rie $f$. L'application ainsi d\Ž finie de
$A[[Y]]$ dans elle-m\ me, $a
\mapsto f.a$, est lin\Ž aire. D\Ž signons-la
par $\rD_f$ (sans grand risque de confusion).

\

\su{7.1.1. C'est une d\Ž rivation de l'alg\ bre de Lie
$A[[Y]]$} $<8>$.

\

\noi \bf En effet\rm, il faut montrer que l'on a
$$f.[a,b] = [f.a,b] + [a,f.b]$$
pour tous $a$ et $b$ dans $A[[Y]]$. Il suffit de le faire
lorsque $a$, $b$ et $f$ sont des polyn\™ mes homog\ nes.

\

Soient donc
$$a \in A[Y]_m \ , \ b \in A[Y]_n \ , \ f \in S_r(Y),$$
$$u \in L_m(Y;A) \ , \ v \in L_n(Y;A).$$
Supposons que $a$ est d\Ž termin\Ž \ par $u$, et $b$ par $v$.
Alors $[a,b]$ est d\Ž termin\Ž \ par $w = [u,v]$ (voir au \bf
5.2\rm). Posons $c = f.[a,b]$. On a $[a,b] \in A[Y]_{m+n}$
et $c \in A[Y]_{m+n+r-1}$. De plus,
$$c(y) = \sum_i [y,\dots,f(y)\urcorner^i,\dots,y),b(y)] +
\sum_i [a(y),v(y,\dots,f(y)\urcorner^i,\dots,y)] =$$
$$[f.a,b](y) + [a,f.b](y).\qed$$
L'application $f \mapsto \rD_f$ de $S(Y)$ dans $\Der(A : Y)$
est lin\Ž aire.

\

\su{7.1.2. C'est un homomorphisme d'alg\ bres de Lie} $<9>$.

\

\noi \bf En effet\rm, il faut montrer que l'on a
$$\rD_{[f,g]} = [\rD_f,\rD_g],$$
autrement dit, que l'on a
$$[f,g].a = f.(g.a) - g.(f.a)$$
pour tous
$$a \in A[[Y]] \ , \ f \in S(Y) \ , \ g \in S(Y).$$
Or, cela d\Ž coule du lemme \bf 4.2\rm.\qed

\

\su{7.2. Homomorphisme de $B$ dans $\Der(A : Y)$}

\

Reprenons l'action formelle $\rd : B \to S(Y)$. Pour chaque $b
\in B$, on a $\rd_b \in S(Y)$. Nous venons de voir comment
$\rd_b$ d\Ž finit une d\Ž rivation de l'alg\ bre de Lie
$A[[Y]]$. D\Ž signons cette d\Ž rivation par $\sigma(b)$.

\

L'application compos\Ž e $\sigma : B \to \Der(A : Y)$ est donc
un homomorphisme d'alg\ bres de Lie, compos\Ž \ de l'action
formelle $\rd : B \to S(Y)$ suivie de l'homomorphisme $S(Y)
\to
\Der(A : Y)$ du \bf 7.1.2\rm.

\

\su{7.3. Crochet sur $W = A[[Y]] \times B$}

\

Etant donn\Ž s deux \Ž l\Ž ments $(f,b)$ et $(g,c)$ de $W$,
on pose
$$[(f,b),(g,c)] = ([f,g] + \rd_b . g - \rd_c .
f,[b,c]).$$ Cela d\Ž finit sur $W$ une structure d'alg\ bre
de Lie qui n'est autre que le \bf produit semi-direct \rm de
l'alg\ bre de Lie $B$ par l'alg\ bre de Lie $A[[Y]]$
relativement \ˆ \ l'homomorphisme
$\sigma : B \to \Der(A : Y)$ introduit ci-dessus.

\

\

\

\head 8. Action triangulaire\endhead

\

\eightpoint

On se donne  

$\rD$ : une action de $A$ sur $X$,

$\rd$ : une action de $B$ sur $Y$,

$Z = X \times Y$ : l'espace vectoriel produit,

$W = W(A,B;\rd)$ : le produit d'entrelacement.

On  fait agir $W$ sur $Z$ canoniquement, \sl en cascade\rm.
C'est cette action $\Delta$ que l'on baptisera \bf action
triangulaire\rm. C'est, en quelque sorte, un produit
d'entrelacement $\Delta = \Delta(\rD,\rd)$ de l'action $\rd :
B
\to S(Y)$ par l'action
$\rD : A \to S(X)$.

\tenpoint

\

Comme on l'a vu ci-dessus (au \bf 5.6)\rm, l'action $\rD : A
\to S(X)$ se prolonge en une action $\rD : A[[Y]] \to S(Z)$.

\

D'autre part l'alg\ bre de Lie $S(Y)$ est plong\Ž e
canoniquement dans l'alg\ bre de Lie $S(Z)$ (voir ci-dessus
au \bf 5.7.2\rm) donc
$\rd : B \to S(Y) \to S(Z)$ est une action de
$B$ sur
$Z$.

\

Pour chaque couple $(f,b) \in A[[Y]] \times B$, on posera
$$\Delta_{(f,b)} = \rD_f + \rd_b \ , \ \text{c'est un \Ž l\Ž
ment de} \ \ S(Z).$$ 
On obtient ainsi une application $\Delta : W \to S(Z)$ qui
est lin\Ž aire.

\

On a alors le r\Ž sultat tr\ s important suivant.

\

\su{8.1. Th\Ž or\ me} $<12>$.

\

\sl L'application $\Delta : W \to S(Z)$
est une action formelle du produit d'entrelacement
$W(A,B;\rd)$ sur l'espace produit $Z = X \times Y$\rm.

\

\su{D\Ž monstration} Le seul point d\Ž licat consiste \ˆ \ v\Ž
rifier que l'on bien
$$[\Delta_{(f,b)},\Delta_{(g,c)}] =
\Delta_{[(f,b),(g,c)]}.$$
Le premier membre s'\Ž crit
$$[\rD_f + \rd_b,\rD_{g} + \rd_c] = [\rD_f,\rD_g] +
[\rD_f,\rd_c] + [\rd_b,\rD_g] +  [\rd_b,\rd_c].$$ 
D'autre part, on a
$$[(f,b),(g,c)] = ([f,g] + \rd_b g - \rd_{c} f,[b,c]).$$
De sorte que le second membre s'\Ž crit
$$\Delta_{[(f,b),(g,c)]} = \rD_{[f,g]} + \rD_{\rd_b g} -
\rD_{d_c f} + \rd_{[b,c]}.$$
Il suffira donc de v\Ž rifier, successivement, que l'on a
$$\rD_{[f,g]} = [\rD_f,\rD_g].\tag 1$$
$$\rd_{[b,c]} = [\rd_b,\rd_c].\tag 2$$
$$\rD_{\rd_b f} - \rD_{\rd_c g} = [\rd_b,\rD_g] -
[\rd_c,\rD_f].\tag 3$$
Or, le (1) d\Ž coule du \bf 5.6 \rm ci-dessus. Le (2) d\Ž
coule du fait que $\rd$ est une action formelle. Quant au
(3), on montrera d'abord ceci : pour \bf tout \rm couple
$(a,b) \in W$, on a
$$\rD_a.\rd_b = 0,\tag 4$$
$$\rd_b.\rD_a = \rD_c \ \ \text{o\ } \ \ c = \rd_b.a =
\sigma(b).a.\tag 5$$
\bf En effet \rm : (4) La s\Ž rie formelle $\rd_b$ ne d\Ž
pend pas de $X$ et la s\Ž rie formelle $\rD_a$ prend ses
valeurs dans $X$. Donc $\rD_a.\rd_b = 0$,
(voir au \bf 3.7 \rm ci-dessus).

\

(5) On se sert de la remarque ci-dessus (au \bf 3.5\rm).
Calculer la d\Ž riv\Ž e $\rd_b.\rD_a$ dans $S(Z)$ o\ \
$$\rd_b \in S(Y) \ \ \text{et} \ \ \rD_a \in S(X)[[Y]],$$
revient \ˆ \ la calculer comme d\Ž riv\Ž e sur $Y$. Or, c'est
la d\Ž riv\Ž e de la s\Ž rie \sl compos\Ž e \rm de $a$ suivi
de l'application \bf lin\Ž aire \rm $\rD$, de sorte que (voir
ci-dessus au \bf 3.8) \rm
$\rd_b.\rD_a =
\rD_{\rd_b.a} = 
\rD_c$, comme annonc\Ž \ . \qed

\

On a ainsi
$$[\rd_b,\rD_a] = \rd_b.\rD_a - \rD_a.\rd_b = \rD_b.\rD_a =
\rD_{\rd_b.a}.\tag 6$$
De (6) d\Ž coule imm\Ž diatement (3), ce qui ach\ ve la d\Ž
monstration.\qed

\

\su{8.2. L'action triangulaire} Reprenons la formule qui
d\Ž finit cette action du produit d'entrelacement 
$W = W(A,B;\rd)$ sur l'espace vectoriel produit$Z = X \times
Y$ :
$$\Delta_{(a,b)} (x,y) = \rD_{a(y)}(x) + \rd_b(y).$$
Cette action est \bf triangulaire \rm dans le sens o\ \ elle
comporte trois temps :
l'action $\rd$, en position $b$, commence par agir sur le
point
$y$ de
$Y$ pour donner $\rd_b(y)$ puis l'\Ž l\Ž ment
$a$ de
$A[[Y]]$ agit sur le point
$y$ de $Y$ pour fournir $a(y)$ ce qui enclenche l'action $\rD$
en position $a(y)$ et fait agir $\rD_{a(y)}$ sur le point $x$
de
$X$ pour donner $\rD_{a(y)}(x)$.

\

\noi De mani\ re imag\Ž e, on
peut dire que c'est une action
\sl en cascade\rm.

\

\

\centerline{\bf A partir d'ici on suppose que le corps $K$ est
de caract\Ž ristique nulle\rm.}

\

\

\head 9. Action fondamentale d'une alg\ bre de Lie sur
elle-m\ me\endhead

\

\eightpoint

Pour le crochet $[\xi,\eta] = \xi.\eta - \eta.\xi$, on sait
d\Ž j\ˆ \  que $S(B)$ est une alg\ bre de Lie (voir au \bf 4
\rm ci-dessus). On introduit, \ˆ \ pr\Ž sent, l'\bf action
fondamentale \rm de $B$ sur elle-m\ me : c'est une action
formelle particuli\ re, un homomorphisme d'alg\ bre de Lie
particulier
$\rd : B
\to S(B)$. 
\tenpoint

\

C'est \ˆ \ dessein que l'on choisit de d\Ž finir cette action
fondamentale pour l'alg\ bre
$B$ plut\™ t que $A$, afin de faciliter la transition entre
ce paragraphe \bf 9 \rm et le paragraphe \bf 10\rm, suivant.

\

On proc\Ž dera par \Ž tapes, comme suit.

\

\su{9.1. Une s\Ž rie g\Ž n\Ž ratrice}

\

On commence par consid\Ž rer la \lq\lq s\Ž rie g\Ž n\Ž
ratrice" suivante :
$$G(T) = \frac {Te^T}{e^T - 1} = \sum_n t_n T^n.$$
Les coefficients $t_n$ sont des nombres
rationnels qui appartiennent donc au corps
$K$ qui est de caract\Ž ristique nulle ! 

\

Plus pr\Ž cis\Ž
ment, on a
$$t_0 = 1 \ , \ t_1 = \frac{1}{2} \ \ \text{et, pour} \ \ n
\geq 1
\ , \ t_{2n} = \frac{b_{2n}}{(2n)!} \ , \ t_{2n+1} = 0,$$
o\ \ les $b_{2n}$ sont les nombres de \smc BERNOULLI\rm.
Attention, cependant, on a $t_1 = - b_1$.

\

\su{9.2. Convention} Comme d'habitude, on prolonge la suite
des coefficients $t_n$ par la convention suivante : $t_i = 0$
pour $i < 0$.

\

\su{9.3. D\Ž finition de l'action fondamentale}

\

Pour chaque $b \in B$ et chaque $n$, on consid\ re
l'application $u_{b,n} : B^n \to B$
d\Ž finie par
$$u_{b,n}(y_1,\dots,y_n) = t_n (\rad \ y_1) \circ (\rad \
y_2)
\circ \dots \circ (\rad \ y_n)(b).$$
Ici, comme d'habitude,
$$\rad \ y : B \to B$$
d\Ž signe l'application lin\Ž aire \bf adjointe\rm
$$(\rad \ y)(b) = [y,b].$$
On pose aussi
$$\rd_{b,n}(y) = t_n (\rad \ y)^n(b),$$
de sorte que
$$\rd_{b,n}(y) = u_{b,n}(y,\dots,y).$$
Ainsi $\rd_{b,n}$ est un polyn\™ me homog\ ne de
degr\Ž
\
$n$, \ˆ \ variables et coefficients dans $B$, d\Ž termin\Ž \
par $u_{b,n} \in L_n(B;B)$. 

\

\su{9.3.1} Autrement dit, on a $\rd_{b,n} \in S(B)_n$. $<13>$.

\

On d\Ž signe enfin par
$$\rd_b = \sum_n \rd_{b,n}$$
la s\Ž rie formelle correspondante, qui appartient \ˆ \
$S(B)$.

\

On d\Ž finit ainsi une application \bf lin\Ž aire\rm,
canonique,
$$\rd : B \to S(B) \ , \ b \mapsto \rd_b,$$
que l'on appellera \bf action fondamentale \rm de $B$.

\

Le th\Ž or\ me important suivant montre que cette action
est une action formelle au sens donn\Ž \ ci-dessus, au \bf
5\rm.

\

\su{9.4. Th\Ž or\ me} $<14>$.

\

\sl Pour toute alg\ bre de Lie $B$, son action fondamentale
$\rd : B \to S(B)$ est un homomorphisme d'alg\ bres de
Lie\rm.

\

Autrement dit, $\rd$ est une action formelle de $B$ sur
elle-m\ me.

\

Nous d\Ž montrons ce th\Ž or\ me plus bas (au \bf 9.10\rm).
Il nous faut d'abord \Ž tablir un certain nombre de r\Ž
sultats auxiliaires afin de faciliter la d\Ž monstration
finale.

\

\noi On commence par v\Ž rifier l'indentit\Ž \ remarquable
suivante pour la fonction g\Ž n\Ž ratrice $G$.

\

\su{9.5. Lemme} La s\Ž rie formelle $G$ v\Ž rifie l'identit\Ž
\ suivante :
$$G(x + y) = L(x,y) + L(y,x)$$
o\ 
$$L(x,y) = \frac{G(x + y) - G(y)}{x} G(x).$$
\bf En effet\rm, on a
$$G(x) = \frac{xe^x}{e^x - 1} \ , \ G(x+y) =
\frac{(x+y)e^{x+y}}{e^{x+y} -1} \ , \ G(x)G(y) =
\frac{xye^{x+y}}{(e^x - 1)(e^y - 1)}.$$
Ainsi
$$U := \frac{G(x)}{x} + \frac{G(y)}{y} = \frac{e^x}{e^x-1} +
\frac{e^y}{e^y-1} = \frac{2 e^{x+y} - e^x -
e^y}{(e^x-1)(e^y-1)}$$
et
$$V := \frac{G(x)G(y)}{x} + \frac{G(x)G(y)}{y} =
\frac{(x+y)e^{x+y}}{(e^x-1)(e^y-1)}.$$
De sorte que
$$L(x,y) + L(y,x) = U G(x+y) - V =
\frac{(x+y)e^{x+y}}{(e^x-1)(e^y-1)}.\left(\frac{2 e^{x+y} -
e^x - e^y}{e^{x+y} -1} - 1\right) =$$
$$\frac{(x+y)e^{x+y}}{e^{x+y} -1}.\frac{e^{x+y} - e^x -
e^y + 1}{(e^x-1)(e^y-1)} = G(x + y).\qed$$

\

A pr\Ž sent, on va \Ž tablir deux autres identit\Ž s. Elles
sont
\sl combinatoires \rm et font intervenir les coefficients
$t_n$. La premi\ re est tr\ s simple et doit nous servir
plus loin.

\

\su{9.6. Lemme} Pour tout entier naturel $m$, on a
$$\frac{1}{m!} = \sum_{n + r = m}
\frac{t_r}{(n+1)!}.$$
\bf En effet\rm, on a
$$e^x = \frac{e^x - 1}{x} G(x),$$
donc
$$\sum_m \frac{x^m}{m!} = \left(\sum_n
\frac{x^n}{(n+1)!}\right) \left( \sum t_rx^r\right) = \sum_m
\left(\sum_{n+r = m}
\frac{t_r}{(n+1)!}\right) x^m.$$
Le r\Ž sultat s'obtient par identification.\qed

\

La seconde identit\Ž \ combinatoire n\Ž cessite davantage de
calculs.

\

\su{9.7. Lemme} Pour tous entiers $m,p,q,$ tels que $m = p+q$,
on a
$$\binom mp t_m = \sum \Sb n+r = m+1 \\ p \leq n - 1 \endSb
\binom np t_n t_r + \sum \Sb n + r = m +1 \\ q \leq r - 1
\endSb \binom rq t_n t_r.$$

\

\su{D\Ž monstration} On utilise la m\Ž thode classique des s\Ž
ries g\Ž n\Ž ratrices. 

\

On pose 
$$H(x,y) = \sum_{p+q = m}  \binom mp 
t_m x^p y^q,$$
$$I(x,y) =  \sum_{p+q = m}\left(\sum \Sb
n + r = m+1 \\ p \leq n-1 \endSb \binom np t_n t_r\right) x^p
y^q,$$
$$J(x,y) = \sum_{p+q = m}\left(\sum \Sb
n + r = m+1 \\ q \leq r-1 \endSb \binom rq t_n t_r\right)
x^px^q.$$
Il suffira donc de montrer que l'on a
$$H(x,y) = I(x,y) + J(x,y).$$
\bf Calcul de $H$\rm. On a
$$H(x,y) = \sum_m t_m (x+y)^m = G(x+y).$$
\bf Calcul de $I$\rm. On a
$$yI(x,y) = \sum_{p+q = m} \left(\sum \Sb n + r = m + 1 \\ p
\leq n - 1 \endSb \binom np t_n t_r\right) x^p y^{q+1}.$$
Or, pour $n+r=m+1$ et $p+q = m$, on a $q+1 = m+1-p = n+r-p$.
Ainsi,
$$yI(x,y) = \sum \Sb n + r = m + 1 \\ p \leq n-1 \endSb
\binom np t_n t_r x^p y^{n-p} y^r = \sum_{n+r = m+1} t_n t_r
((x+y)^n - x^n)y^r =$$
$$\sum_r (G(x+y) - G(x))t_ry^r = (G(x+y) - G(x))G(y) =
G(x+y)G(y) - G(x)G(y).$$
\bf Calcul de $J$\rm. Par sym\Ž trie, on a
$$J(x,y) = I(y,x).$$
D'o\ \
$$xJ(x,y) = G(x+y)G(x) - G(x)G(y).$$
Ainsi,
$$I(x,y) + J(x,y) = \frac{G(x+y)G(x) - G(x)G(y)}{x} +
\frac{G(x+y)G(y) - G(x)G(y)}{y} = G(x+y)$$
d'apr\ s le lemme \bf 9.5\rm. D'o\ \ l'identit\Ž \ annonc\Ž
e
$$H(x,y) = I(x,y) + J(x,y).\qed$$

\

Voici, \ˆ \ pr\Ž sent, des r\Ž sultats auxiliaires sur les
d\Ž rivations des alg\ bres de Lie. Ils nous servirons \ˆ \
\Ž tablir le th\Ž or\ me \bf 9.4 \rm ainsi que d'autres r\Ž
sultats encore, plus loin.

\

\su{9.8. Lemme} Soit $\rD$ une d\Ž rivation quelconque de
l'alg\ bre de Lie $B$. Alors, pour tous \Ž l\Ž ments $a \in
B \ , \ b \in B$, et tout entier $m \geq 0$, on a les
identit\Ž s suivantes :
$$\sum_{0 \leq k \leq m} \rD^k[a,\rD^{m-k} b] = \sum_{0 \leq i
\leq m} \binom {m+1}{i+1} [\rD^ia,\rD^{m-i}b] = \sum \Sb i +
j = m \\ i \geq 0 \endSb \binom {m+1}{i+1}
[\rD^ia,\rD^jb].\tag 1$$
$$\sum_{0 \leq k \leq m} \left(\rD^k[a,\rD^{m-k} b] +
\rD^k[\rD^{m-k}a, b]\right) = \sum_{0 \leq i \leq m} \binom
{m+2}{i+1} [\rD^ia,\rD^{m-i}b].\tag 2$$

\

\su{D\Ž monstration} 

\

\noi (1) On utilise la formule de Leibniz pour
les d\Ž riv\Ž es d'ordre sup\Ž rieur. Il vient
$$\sum_{0 \leq k \leq m} \rD^k[a,\rD^{m-k} b] = \sum_{0 \leq
k \leq m} \sum_i \binom ki [\rD^ia,\rD^{m-i}b].$$
Or, pour $i \geq 0$ fix\Ž , on a
$$\sum_{0 \leq k \leq m} \binom ki = \binom 0i + \binom 1i +
\dots + \binom mi = \binom {m+1}{i+1}.$$
D'o\ \ la premi\ re identit\Ž .

\

\noi (2) D'apr\ s le (1), en \Ž changeant les r\™ les de $a$
et de $b$, on a
$$\sum_{0 \leq k \leq m} \rD^k[\rD^{m-k}a,b] = \sum \Sb i +
j = m \\ j \geq 0 \endSb \binom {m+1}{j+1}
[\rD^ia,\rD^jb] = \sum_{0 \leq i \leq m} \binom {m+1}{i}
[\rD^i a,\rD^{m-i}b].$$
Donc
$$\sum_{0 \leq k \leq m} \left(\rD^k[a,\rD^{m-k} b] +
\rD^k[\rD^{m-k}a, b]\right) = \sum_{0 \leq i \leq m} \left(
\binom {m+1}{i+1} + \binom {m+1}{i}\right) [\rD^ia,\rD^{m-i}b]
=
$$
$$\sum_{0
\leq i \leq m}  \binom {m+2}{i+1}
[\rD^ia,\rD^{m-i}b].\qed$$

\ 

Combinant ce dernier lemme et l'identit\Ž \ combinatoire du
\bf 9.7\rm, on \Ž tablit le r\Ž sultat suivant.

\

\su{9.9. Lemme} Soit $\rD$ une d\Ž rivation quelconque de
l'alg\ bre de Lie $B$. Alors, pour tous $a \in B$ et $b \in
B$, et tout entier $m \geq 0$, on a
$$t_m \rD^m[a,b] = \sum_{n+r = m+1} t_nt_r \left(\sum_{0
\leq k \leq r-1} \rD^k [\rD^na,\rD^{r-k-1}b] + \sum_{0 \leq
k \leq n-1} \rD^k[\rD^{n-k-1}a,\rD^rb] 
\right).$$

\

\su{D\Ž monstration} La formule de Leibniz donne
$$\rD^m[a,b] = \sum_{p+q = m} \binom mp [\rD^pa,\rD^qb].$$
D'autre part, l'identit\Ž \ \bf 9.8\rm(1) \rm permet d'\Ž
crire
$$\sum_{0 \leq k \leq r-1} \rD^k[\rD^na,\rD^{r-k-1}b] =
\sum_{0 \leq i \leq r-1} \binom r{i+1}
[\rD^{n+i}a,\rD^{r-i-1}b].$$
Le terme en  $[\rD^pa,\rD^qb]$ dans la somme ci-dessus
correspond \ˆ \ 
$$n+i = p \ , \ r-i-1 = q \ , \ q \leq r-1.$$
Son coefficient est donc $\binom rq$. De sorte
que l'on a
$$\sum_{0 \leq k \leq r-1} \rD^k[\rD^na,\rD^{r-k-1}b] =
\sum \Sb p + q = r - 1 \\ q \leq r - 1 \endSb \binom rq
[\rD^pa,\rD^qb].$$  
De m\ me, on a
$$\sum_{0 \leq k \leq n - 1} \rD^k[\rD^{n-k-1}a,\rD^r b] =
\sum \Sb p+q = n-1 \\ p \leq n-1 \endSb \binom np
[\rD^pa,\rD^qb].$$
Ainsi, dans la somme suivante
$$\sum_{n+r = m+1} t_nt_r \left(\sum_{0
\leq k \leq r-1} \rD^k [\rD^na,\rD^{r-k-1}b] + \sum_{0 \leq
k \leq n - 1} \rD^k[\rD^{n-k-1}a,\rD^rb] 
\right),$$
le coefficient de $[\rD^pa,\rD^qb]$ est \Ž gal \ˆ \
$$\sum \Sb n + r = m + 1 \\ q \leq r - 1 \endSb \binom rqt_n
t_r \  + \sum \Sb n + r = m + 1 \\ p \leq n-1 \endSb \binom
np t_n t_r$$
qui est \Ž gal \ˆ \ $\binom mp t_m$,
d'apr\ s le \bf 9.7 \rm ci-dessus, ce qui ach\ ve la d\Ž
monstration.\qed

\

\su{9.10. D\Ž monstration du th\Ž or\ me 9.4}

\

Il faut montrer que, pour tous $a \in B$, $b \in B$, on a
$$\rd_{[a,b]} = [\rd_a,\rd_b].$$
On proc\ de par composantes homog\ nes.

\

Soient $f_m$ et $g_m$, respectivement, les composantes homog\
nes de degr\Ž \ $m$ du premier et du second membre. Il
s'agit de montrer que l'on a
$$f_m(y) = g_m(y) \ \ \text{pour tout} \ \ y \in B.$$
Posons $\rad \ y = \rD$. Il vient
$$f_m(y) = \rd_{[a,b]} (y) = t_m (\rad \ y)^m [a,b] = t_m
\rD^m [a,b].$$
D'autre part, on a
$$g_m(y) = \sum_{n + r = m + 1} [\rd_{a,n},\rd_{b,r}](y) =
\sum_{n + r = m + 1} (\rd_{a,n}.\rd_{b,r} -
\rd_{b,r}.\rd_{a,n})(y).$$
\bf Calcul de la d\Ž riv\Ž e $(\rd_{a,n}.\rd_{b,r})(y)$\rm.
On se sert de l'application multilin\Ž aire $u_{b,r}$ qui d\Ž
termine le polyn\™ me homog\ ne $\rd_{b,r}$ (voir ci-dessus,
au
\bf 9.3\rm)
$$\rd_{b,r}(y) = u_{b,r}(y,\dots,y) = t_r \rD^r b,$$
et de la valeur
$$\rd_{a,n}(y) = t_n \rD^n a.$$
Il vient
$$(\rd_{a,n}.\rd_{b,r})(y) = \sum_{1 \leq k \leq r} u_{b,r}
(y,\dots,t_n \rD^n a\urcorner^k,\dots,y) =  \sum_{1 \leq k
\leq r} t_n  u_{b,r} (y,\dots,\rD^n a\urcorner^k,\dots,y)$$
Or,
$$u_{b,r} (y,\dots,\rD^n a\urcorner^k,\dots,y) = t_r (D \circ
\dots \circ \rD^n a\urcorner^k \circ \dots \circ \rD)(b) =
t_r \rD^{k-1}[\rD^na,\rD^{r-k}b].$$ 
D'o\ \
$$(\rd_{a,n}.\rd_{b,r})(y) = \sum_{1
\leq k \leq r} t_n t_r \rD^{k-1}[\rD^na,\rD^{r-k}b].$$
\bf De m\ me\rm
$$(\rd_{b,r}.\rd_{a,n})(y) = \sum_{1
\leq k \leq n} t_n t_r \rD^{k-1}[\rD^rb,\rD^{n-k}a].$$
\bf Ainsi\rm,
$$g_m(y) = \sum_{n+r = m+1} t_nt_r \left(\sum_{1
\leq k \leq r} \rD^{k-1} [\rD^na,\rD^{r-k}b] + \sum_{1 \leq
k \leq n} \rD^{k-1}[\rD^{n-k}a,\rD^rb] 
\right).$$
On a donc, d'apr\Ž s le \bf 9.9 \rm ci-dessus,
$$g_m(y) = t_m \rD^m[a,b] = f_m(y).\qed$$

\

\su{9.11. Sur l'origine de la notion d'action fondamentale}

\

Soit $B$ une alg\ bre de Lie \bf norm\Ž e
compl\ te\rm, r\Ž elle ou complexe. 

\

On sait lui associer le \bf groupuscule de Lie
d\Ž fini par $B$ \rm (voir Bourbaki [4], III, pages 168-169,
dont on gardera, ici, les
notations, autant que possible, sauf \ˆ \ remplacer le $A$ par
$B$).

\

Soit $G$ ce groupuscule : c'est un voisinage ouvert de $0$
dans $B$. Bien entendu, l'alg\ bre de Lie $L(G)$ du
groupuscule $G$ est identifi\Ž e \ˆ \ $B$.

\

Il existe, par hypoth\ se, un \bf morceau \rm de loi d'op\Ž
ration \ˆ \ droite, analytique, canonique, du groupuscule $G$
sur la vari\Ž t\Ž \ $G$ : c'est une application,
partiellement d\Ž finie, $G \times G \to G$, $(x,y) \mapsto
x.y$.

\

A ce morceau de loi, correspond une \bf loi d'op\Ž ration
infinit\Ž simale \rm \ˆ \ droite, analytique, de $B = L(G)$
dans $G$ (voir ibid. III, page 165, 18.7).

\

A chaque loi infinit\Ž simale nous avons associ\Ž \
(ci-dessus, au paragraphe \bf 6\rm) une action formelle de
l'alg\ bre de Lie $B$ sur l'espace vectoriel $B$.

\

En l'occurence, on a le r\Ž sultat suivant.

\

\su{9.12 Th\Ž or\ me} \sl L'action formelle
associ\Ž e
\ˆ
\ la loi d'op\Ž ration infinit\Ž simale de $B$ dans le
groupuscule
$G$ n'est autre que l'action fondamentale $\rd : B \to S(B)$
que nous venons de d\Ž finir\rm. $<15>$. 

\

\su{On le v\Ž
rifie, succintement, comme suit}

\

Voici quelques d\Ž tails suppl\Ž
mentaires afin de faciliter le
\sl raccordement
\rm de cette face \bf analytique \rm avec l'aspect \bf formel
\rm que nous avons pr\Ž sent\Ž .

\  

Dans le groupuscule $G$, le produit $x.y$ est d\Ž fini par la
s\Ž rie de Hausdorff (voir Bourbaki [4], II, pages 55-57 et
l'exercice 3, page 9, ainsi que III.4.2, page 168) :
$$x.y = H(x,y) = x + y + \frac{1}{2} [x,y] + \frac{1}{12}
[x,[x,y]] + \frac{1}{12}
[y,[y,x]] + \cdots.$$
L'application \sl
exponentielle \rm de l'alg\ bre de Lie $B$ dans son
groupuscule n'est autre que $f = \roman{id}_B : B \to
B$. Que dire de la \sl diff\Ž rentielle
\rm
$$(\rD_af)(b) = \lim \frac{e^{ta}.b - b}{t} = \frac{(ta).b -
b}{t} \ ?$$
On a 
$$(ta).b - b = H(ta,b) - b = ta + \frac{t}{2} [a,b] +
\frac{t^2}{12} [a,[a,b]] + \frac{t}{12}
[b,[b,a]] + \cdots.$$
Autrement dit, $(\rD_af)(b)$ est donn\Ž \ par la s\Ž rie
$H_1(a,b)$,  somme des termes de $H(a,b)$ dont le degr\Ž \ en
$a$ est $1$ : c'est-\ˆ \ -dire
$$(\rD_af)(b) = H_1(a,b)$$
o\ \ (voir Bourbaki [4], II.6,
exercice 3d, page 90),
$$H_1(a,b) = a + \frac{1}{2} [a,b] + \sum_{n \geq 1}
\frac{1}{(2n)!} b_{2n} (\rad \ b)^{2n} (a),$$
les $b_{2n}$ \Ž tant les nombres de BERNOULLI. Le th\Ž
or\ me en d\Ž coule.\qed

\

Avec l'exemple originel (au paragraphe \bf 6\rm), c'est cela
qui a \Ž t\Ž \ le point de d\Ž part de notre d\Ž veloppement.
C'est l\ˆ \
que nous avons puis\Ž \ une part de cette inspiration qui
nous a conduit \ˆ \ la notion de produit d'entrelacement et
d'action formelle.

\

\su{9.13. Encore un petit mot de commentaire}

\

Lorsque $B$ est
une alg\ bre de Lie norm\Ž e compl\ te, r\Ž elle ou
complexe, l'endomorphisme
$(\rad \ y)$ est un op\Ž rateur continu de l'espace $B$.  La
s\Ž rie enti\ re
$\sum_n t_n (\rad \ y)^n$ est alors
normalement convergente, avec un rayon de convergence infini.
Elle d\Ž finit donc une fonction analytique sur $B$,
\ˆ \ valeurs dans l'espace des op\Ž rateurs continus de
l'espace vectoriel $B$. Reprenant la s\Ž rie g\Ž n\Ž ratrice
$G$ du
\bf 9.1\rm, on peut ainsi \Ž crire
$$G(\rad
\ y) = \sum_n t_n (\rad \ y)^n.$$

\

Dans le cas g\Ž n\Ž ral, \bf purement alg\Ž brique\rm, on
peut encore
\Ž crire, \bf formellement\rm,
$$\rd_b(y) = \sum_n t_n (\rad \ y)^n(b) = G(\rad \
y) (b).$$
On va donner une justification de cette \Ž criture et
un mot d'explication au sujet de la nature de la s\Ž
rie 
$$G(\rad \
y) = \sum_n t_n (\rad \ y)^n.$$
Pour $b \in B$ fix\Ž , la s\Ž rie $\sum_n t_n (\rad \ y)^n(b)$
appartient
\ˆ
\ $S(B)$. Comme fonction, \ˆ \ la fois, de $y$ et de $b$,
elle appartient \ˆ
\ $S(B)[[B]]$ lequel est identifi\Ž \ \ˆ \ $B[[B,B]]$. 

\

\noi \{En posant 
$$w(y_1,\dots,y_n,b) = t_n (\rad \ y_1) \circ (\rad \
y_2)
\circ \dots \circ (\rad \ y_n)(b),$$
$$u_b(y_1,\dots,y_n) = v_{(y_1,\dots,y_n)}(b) = t_n (\rad \
y_1) \circ (\rad \ y_2)
\circ \dots \circ (\rad \ y_n)(b),$$
on obtient une application multilin\Ž aire $w \in
L_{(n,1)}(B,B;B)$ qui s'identifie, d'une part \ˆ \ 
l'application $u
\in L_1(B;L_n(B;B))$ et, d'autre part, \ˆ \ $v \in
L_n(B;L_1(B;B))$. De sorte que $w(y,\dots,y,b)$ est un
polyn\™ me homog\ ne de bidegr\Ž \ $(n,1)$ et appartient
ainsi \ˆ \ $B[B,B]_{(n,1)}$, tandis que les
polyn\™ mes homog\ nes $u_b(y,\dots,y)$ et
$v_{(y,\dots,y)}(b)$ appartiennent, respectivement, \ˆ \
$S(B)_n[B]_1$ et $S(B)_1[B]_n$.\}

\

\noi Pour  $y \in B$ donn\Ž , $(\rad \ y)$ est
un endomorphisme de l'espace
vectoriel
$B$, autrement dit, un op\Ž rateur lin\Ž aire
$(\rad \ y) \in L_1(B,B) = \roman{End}(B)$. Ainsi, $(\rad \
y)$, comme fonction de
$y$, est un polyn\™ me homog\ ne de degr\Ž \ $1$ \ˆ \
variables dans $B$ et coefficients dans $\roman{End}(B)$;
c'est un \Ž l\Ž ment de $\roman{End}(B)[B]_1$. Plus g\Ž n\Ž
ralement, on a $(\rad \ y)^n \in \roman{End}(B)[B]_n$ donc
$\sum_n t_n (\rad \ y)^n \in \roman{End}(B)[[B]] \subset
S(B)[[B]]$.

\

\noi Lorsque  
l'op\Ž rateur
$(\rad
\ y)$ est
\bf nilpotent\rm,  $\sum_n t_n (\rad \ y)^n$ est une somme
\bf finie \rm et repr\Ž sente un \Ž l\Ž ment de
$\roman{End}(B)$ que l'on peut encore d\Ž signer, sans grand
danger,  par
$G(\rad
\ y)$. En particulier, lorsque l'alg\ bre de Lie $B$ est
\bf nilpotente\rm,  $G(\rad \ y)$ est un polyn\™ me en $y$
appartenant \ˆ \  $\roman{End}(B)[B]$.

\

\noi Sinon, dans le cas g\Ž n\Ž ral, voici comment justifier
l'\Ž criture 
$$G(\rad
\ y) = \sum_n t_n (\rad \ y)^n.$$
On commence par observer, bri\ vement, ceci. Lorsque $E$ est
une $K$-alg\ bre, o\ \  $\circ$ d\Ž signe la multiplication,
l'ensemble de s\Ž ries formelles $E[[X]]$ h\Ž rite de la
structure de $K$-alg\ bre de $E$, avec la multiplication
naturelle des polyn\™ mes homog\ nes : $(f_n \circ g_r)(x) =
f_n(x) \circ g_r(x)$. Si, de plus, la $K$-alg\ bre $E$
est associative et poss\ de une unit\Ž \ $1$, alors l'anneau
des s\Ž ries formelles $K[[T]]$ est (isomorphe \ˆ ) une
sous-alg\ bre de
$E[[X]]$. La s\Ž rie formelle $e^T$, ainsi que la s\Ž rie g\Ž
n\Ž ratrice
$$G(T) = \frac {Te^T}{e^T - 1} = \sum_n t_n T^n,$$
sont donc des \Ž l\Ž ments de $E[[X]]$.

\

\noi Dans le cas particulier o\ \ $E = X = \roman{End}(B)$,
il devient clair que l'on a $G(T) \in S(\roman{End}(B))$.
Ainsi, $G(\rad \ y)$ n'est autre que la s\Ž rie que l'on
obtient en substituant le polyn\™ me homog\ ne, de
degr\Ž \ $1$, $(\rad
\ y)$
\ˆ
\
$T$ dans la s\Ž rie formelle $G(T)$. On a aussi
$$G(\rad \ y) = \frac{(\rad \ y) e^{(\rad \ y)}}{e^{(\rad \
y)} - 1}.$$
L'exponentielle $e^{\rad \ y}$ s'\Ž
crivant, traditionnellement, $\roman{Ad} \ y$, il vient aussi
$$G(\rad \ y) = \frac{(\rad \ y) \roman{Ad} \ y}{\roman{Ad}
\ y - 1}.$$

\

\

\head 10. Le produit d'entrelacement de deux alg\ bres de
Lie\endhead

\

\eightpoint \`A chaque action formelle donn\Ž e, $\rd : B \to
S(Y)$,  correspond un produit d'entrelacement
$W(A,B;\rd) = A[[Y]] \times B$, comme on l'a dit plus
haut (au paragraphe \bf 7\rm). Parmi tous ces produits
d'entrelacement, il en est un, particulier, que l'on  d\Ž
signera, simplement, par
$W(A,B)$ : c'est celui qui
correspond \ˆ \ l'action fondamentale $\rd : B \to S(B)$ de
l'alg\ bre de Lie
$B$ sur elle-m\ me. D'une certaine mani\ re, il est
\bf intrins\ que\rm. Comme dans le cas des groupes
abstraits, ce produit $W(A,B) = A[[B]] \times B$ agit \sl
en cascade \rm sur l'espace vectoriel produit
$A \times B$, en une action \sl triangulaire\rm.

\tenpoint

\

Lorsque \rm $\rd : B
\to S(B)$ est l'\bf action fondamentale \rm  que l'on vient
de d\Ž finir, on obtient un cas particulier important de
produit d'entrelacement
$W(A,B;d)$ : on appellera ce cas particulier le \bf produit
d'entrelacement (fondamental) \rm  de l'alg\ bre de Lie $B$
par l'alg\ bre de Lie $A$ et on le d\Ž signera, simplement,
par
$W(A,B)$. C'est lui que l'on essayait de d\Ž finir, au d\Ž
part, et lui qui a donn\Ž \ lieu aux pr\Ž sents
d\Ž veloppements.

\

\su{10.1. L'essentiel} Rappelons comment
est construit $W(A,B)$, en se r\Ž
f\Ž rant au paragraphe
\bf 7
\rm ci-dessus.

\

\noi Ici, $\Der(A : B)$ est l'alg\ bre de Lie $\frak d
(A[[B]])$ des d\Ž rivations de l'alg\ bre de Lie $A[[B]]$.

\

\noi A chaque $b \in B$, l'homomorphisme $\sigma : B \to
\Der(A : B)$ associe la d\Ž rivation
$\sigma(b) \in \frak d (A[[B]]$, d\Ž rivation suivant la s\Ž
rie formelle $\rd_b \in S(B)$ :
$$\rd_b(y) = \sum_n t_n (\rad \ y)^n(b) = G(\rad \ y) (b).$$ 

\

\noi Enfin, $W(A,B) = A[[B]] \times B$ est le \bf produit
semi-direct
\rm de l'alg\ bre de Lie $B$ par l'alg\ bre de Lie
$A[[B]]$ relativement \ˆ \  l'homomorphisme $\sigma$. Le
crochet y est d\Ž fini, comme au \bf
7.3 \rm ci-dessus, par la formule suivante : \Ž tant donn\Ž s
deux \Ž l\Ž ments $(f,b)$ et $(g,c)$ de
$W(A,B)$, on a
$$[(f,b),(g,c)] = ([f,g] + \rd_b . g - \rd_c . f,[b,c]),$$
o\ \ le crochet $[f,g]$ est celui que $A[[B]]$ h\Ž rite de
$A$ (voir ci-dessus au \bf 5.1\rm). 

\

\su{10.2. Les d\Ž tails} Afin d'abr\Ž ger, pour d\Ž signer la
d\Ž riv\Ž e d'une s\Ž rie formelle donn\Ž e $f \in
A[[B]]$ suivant la s\Ž rie formelle
$\rd_b$, celle de l'action fondamentale, on \Ž crira $b\star
f$  au lieu de
$\rd_b.f$ :
$$\boxed{b\star f =  \rd_b.f.}$$

\

\noi Toute s\Ž rie formelle $f \in A[[B]]$ se met sous une
forme canonique
$f(x) = \sum_m u_m(x,\dots,x)$ o\ \ $u_m$ est
une application $m$-lin\Ž aire \bf sym\Ž trique\rm, $u_m \in
L_m(A;B)$, (voir ci-dessus au \bf 2.8\rm). Pour la d\Ž riv\Ž
e  de cette s\Ž rie $f$ suivant la s\Ž rie $\rd_b$, on a la
formule suivante (voir ci-dessus au \bf 3.3\rm)
$$(b\star f)(x) = (\rd_b.f)(x) = \sum_m
mu_m(\rd_b(x),x,\dots,x).$$ 
Autrement dit,
$$(b\star f)(x) = \sum_m\sum_n mu_m(t_n(\rad \
x)^n(b),x,\dots,x).$$
On a ainsi le formulaire suivant
$$(b\star f)_m(x) = \sum_{n + r = m + 1} t_nru_r((\rad
\ x)^nb,x,\dots,x).$$
$$[(f,b),(g,c)] = ([f,g] + b\star g - c\star. f,[b,c]),$$
$$[f,g](x) = [f(x),g(x)] \ \ \text{(le crochet de} \ \ A) \ ,
\ [b,c] \ \ \text{(le crochet de} \ \ B).$$

\

\su{10.3. L'action triangulaire} En introduisant \Ž galement
l'action fondamentale
$D : A
\to S(A)$, voici la description de l'action triangulaire 
fondamentale de $W(A,B)$ sur $A \times B$.  

\

\noi C'est, en
quelque sorte, le produit d'entrelacement $\Delta =
\Delta(\rD,\rd)$ des deux actions fondamentales
$\rd : B
\to S(B)$ et $\rD : A \to S(A)$ (voir, ci-dessus,
au paragraphe \bf 8\rm). Cette application $\Delta : W
\to S(A
\times B)$ associe, 
\ˆ
\ chaque couple $(a,b) \in A[[B]] \times B$, l'\Ž l\Ž ment
suivant de  $S(A \times B)$ :
$$\Delta_{(a,b)} = \rD_a + \rd_b.$$
Elle agit de la mani\ re suivante : pour chaque couple
$(x,y) \in A \times B$, on a
$$\Delta_{(a,b)} (x,y) = \rD_{a(y)}(x) + \rd_b(y).$$
Disons, de nouveau, qu'elle est \bf triangulaire \rm dans le
sens o\
\ elle comporte trois temps :
l'action $\rd$, en position $b$, commence par agir sur le
point
$y$ de
$B$ pour donner $\rd_b(y)$ puis l'\Ž l\Ž ment
$a$ de
$A[[B]]$ agit sur le point
$y$ de $B$ pour fournir $a(y) = a_0(y) + a_1(y) + a_2(y) +
\cdots$ ce qui enclenche l'action
$\rD$ en position $a(y)$ et fait agir $\rD_{a(y)}$ sur le
point $x$ de
$A$ pour donner $\rD_{a(y)}(x) = \rD_{a_0(y)}(x) +
\rD_{a_1(y)}(x) + \rD_{a_2(y)}(x) + \cdots\cdot$

\

\noi C'est une action
\sl en cascade\rm, pour ainsi dire.

\

\head Exemples de produits $W(A,B)$\endhead

\

\su{10.4. Un premier exemple} Dans le cas o\ \ les deux
alg\ bres de Lie $A$ et $B$ ont une dimension \Ž gale \ˆ \
$1$,
$A[[B]]$ est idenfi\Ž \ \ˆ \ l'anneau
classique $K[[x]]$ de s\Ž ries formelles, (voir ci-dessus
au \bf 1.3.1\rm), de sorte que l'on a $W(A,B) = K[[x]]
\times K$. L'alg\ bre de Lie $B$ \Ž tant commutative, $(\rad
\ y)$ est nul et l'action fondamentale $\rd_b$ se r\Ž duit
\ˆ \ $\rd_b = b$. La s\Ž rie formelle $f(x) = \sum a_m x^m
\in A[[B]] = K[[x]]$ \Ž tant donn\Ž e, soit $f'(x) = \sum
ma_mx^{m-1}$ sa d\Ž riv\Ž e formelle. On a $(b\star f)(x)
=(d_b.f)(x) = bf'(x)$. Les crochets de $A$ et de $B$ \Ž tant
nuls, tous deux, celui de
$W(A,B)$ est alors donn\Ž
\ par
$$[(f,b),(g,c)] = (b\star g - c\star f,0) = (bg' -
cf',0).$$

\

\su{10.5. Un deuxi\ me exemple} Plus g\Ž n\Ž ralement,
supposons que $B$ soit une alg\ bre de Lie commutative, de
dimension finie $r$, et $A$ une alg\ bre de Lie
de dimension
$1$. Alors $A[[B]]$ est identifi\Ž \ \ˆ \
l'anneau classique de s\Ž ries formelles
$K[[x_1,\dots,x_r]]$. De nouveau, l'action fondamentale
est r\Ž duite \ˆ \ $\rd_b = b$. Pour 
$$b = (b_1,\dots,b_r)
\in B \ , \ x = (x_1,\dots,x_r) \ , \ f(x) \in
K[[x_1,\dots,x_r]],$$
on a
$$(b\star f)(x) = (\rd_b.f)(x) = b_1f'_{x_1}(x) + \dots +
b_rf'_{x_r}(x).$$
Les crochets de $A$ et de $B$ \Ž tant nuls, de
nouveau, celui de $W(A,B)$ s'\Ž crit encore
$$[(f,b),(g,c)] = (b\star g - c\star f,0) =
\left(\sum_i(b_ig'_{x_i} - c_if'_{x_i}),0\right).$$

\

\su{10.6. Un troisi\ me exemple} On suppose que l'alg\ bre
de Lie $B$ est nilpotente et, par exemple, que
$(\rad
\ y)^3$ est nul pour tout $y \in B$. L'action fondamentale
est alors donn\Ž e par 
$$\rd_b(y) = b + \frac{1}{2} [y,b] + \frac{1}{12}
[y,[y,b]].$$
Pour $f = \sum_m u_m(x,\dots,x)$, sous forme canonique, on a
ainsi
$$(b\star f)_0(x) = u_1(b) \ , \ (b\star f)_1(x) =
 2u_2(b,x) + \frac{1}{2}u_1([x,b]),$$
$$(b\star f)_2(x) = 3u_3(b,x,x)
 +  u_2([x,b],x) + \frac{1}{12}u_1([x,[x,b]]),$$
$$(b\star f)_m(x) = \sum_{n + r = m + 1 , n \leq 2}
t_nru_r((\rad
\ x)^nb,x,\dots,x).$$ 
Le crochet de $W(A,B)$ s'\Ž crit 
$$[(f,b),(g,c)] = ([f,g] + b\star g - c\star f,[b,c]).$$

\

\head 11. Repr\Ž sentation des extensions dans le produit
d'entrelacement\endhead

\

\eightpoint Comme dans le cas des groupes abstraits, le r\Ž
sultat suivant illustre la singularit\Ž \ du produit
d'entre\-lacement fondamental $W(A,B)$. Toute alg\ bre de
Lie $C$ qui est une extension de $B$ par $A$ est isomorphe
\ˆ \ une sous-alg\ bre de Lie de $W(A,B)$.

\tenpoint

\

Voici, \ˆ \ pr\Ž sent, un th\Ž or\ me qui montre que l'on
peut repr\Ž senter toute extension $C$ de l'alg\ bre de Lie
$B$ par l'alg\ bre de Lie $A$ dans leur produit
d'entrelacement
$W(A,B)$.

\

Il s'agit de l'analogue pour les alg\ bres de Lie du
premier th\Ž or\ me de Kaloujnine-Krasner pour les groupes
abstraits (voir [6]).

\

Soit $A \to C \overset{\rp} \to \to B$ une extension de
l'alg\ bre de Lie $B$ par l'alg\ bre de Lie $A$. Autrement
dit,
$\rp$ est un homomorphisme surjectif de l'alg\ bre de Lie
$C$ sur l'alg\ bre de Lie
$B$ dont le noyau est la sous-alg\ bre $A = \{c \in
C :
\rp(c) = 0 \}$.

\

On se fixe une application $K$-lin\Ž aire $\rs : B \to C$
quelconque telle que $\rp \circ \rs = \roman{id}_B$.
Autrement dit, $\rs$ est une \bf section lin\Ž aire \rm de
$\rp$. 

\

On va associer \ˆ \ $\rs$ un homomorphisme
$f_{\rs} : C \to W(A,B)$ d'alg\ bres de Lie que l'on
appellera la
\bf repr\Ž sentation associ\Ž e \rm \ˆ \ $\rs$.

\

On pr\Ž sentera cet homomorphisme $f_{\rs}$ par \Ž tapes,
comme suit.

\

\su{11.0. Nota} Afin de rendre plus claire la lecture des
calculs compliqu\Ž s qui vont suivre, on omettra les
parenth\ ses et le symbole $\circ$ de composition des
fonctions partout o\ \ le risque de confusion est \bf
minime\rm. De plus,  pour simplifier, on introduira les
notations suivantes. Pour
$y \in B$ fix\Ž , on pose 
$$D = \rad \ y \ , \ z = \rs
(y) \ , \ R = \rad \ z = \rad \ \rs y.$$
On posera aussi $e = \rs
\rp$, (c'est une application lin\Ž aire $e : C \to C$), et on
observera ceci :
$$\rp e = \rp \ \ \text{et} \ \ e\rs = \rs \ \
\text{puisque} \ \ \rp \rs = \roman{id}_B.$$ 
Ainsi, $D$ est une d\Ž rivation de l'alg\ bre de Lie
$B$ tandis que
$R$ est une d\Ž rivation de l'alg\ bre de Lie $C$.
En particulier, on a $\rp R = D\rp$. On a ainsi, par r\Ž
currence,
$\rp R^n = D^n \rp \ \ \text{pour tout} \ \ n$. On a donc,
plus g\Ž n\Ž ralement, pour tout $k \geq 0$,
$$\rp R^n e R^k = D^n \rp e R^k = D^n \rp R^k = D^n D^k \rp
= D^{n+k} \rp = \rp R^{n+k}.$$ Autrement dit,
$$\rp R^n e R^k = D^{n+k} \rp = \rp R^{n+k}.$$ On se servira
de cette derni\ re identit\Ž , dans les deux sens, comme
d'une
\sl fermeture \Ž clair\rm, de gauche \ˆ \ droite et de droite
\ˆ \ gauche.

\

\

\

\su{11.1. D\Ž finition de $h_c$}

\

Pour $c \in C , \ y \in B$ et $m$ donn\Ž s, on pose 
$$h_{c,m}(y) = \frac{1}{m!} (\rad \ z)^m(c) - \sum_{n+r =
m} \frac{t_r}{(n+1)!} (\rad \ z )^n(\rs \circ \rp) (\rad \
z)^r(c)$$
o\ \ les coefficients $t_r$ sont d\Ž finis par la s\Ž rie
g\Ž n\Ž ratrice $G(T)$ introduite, ci-dessus, au \bf 9.1\rm.
Autrement dit,
$$h_{c,m}(y) = \frac{1}{m!} R^m c - \sum_{n+r =
m} \frac{t_r}{(n+1)!} R^n e R^r c.$$

L'\Ž l\Ž ment $h_{c,m}(y)$ ainsi d\Ž fini appartient
visiblement \ˆ
\ l'alg\ bre de Lie $C$. En fait, il appartient plus pr\Ž
cis\Ž ment
\ˆ \ $A$.

\

\su{11.2. Lemme} $<16>$. Pour tous $c \in C , \ y \in B$,
on a
$$h_{c,m}(y)
\in A.$$
\bf En effet\rm, on calcule 
$$\rp h_{c,m}(y) = \frac{1}{m!} D^m \rp c -
\sum_{n+r=m} \frac{t_r}{(n+1)!} D^n(\rp \rs
\rp)R^r c =$$
$$\frac{1}{m!} D^m \rp c -
\sum_{n+r=m} \frac{t_r}{(n+1)!} D^{n+r} \rp c
= \left(\frac{1}{m!} - \sum_{n+r=m}
\frac{t_r}{(n+1)!}\right) D^m \rp c = 0$$
d'apr\ s le lemme \bf 9.6\rm.\qed

\

Cela prouve que $h_{c,m}$, comme polyn\™ me homog\ ne en
$y$, appartient \ˆ \
$A[B]_m$.

\

On consid\ re la s\Ž rie formelle $h_c = \sum_m h_{c,m}$, de
sorte que l'on a $h_c \in A[[B]]$. On pose enfin
$$f_\rs (c) = (h_c,\rp c).$$
On a ainsi $f_\rs(c) \in A[[B]] \times B = W(A,B)$ et
l'application $f_\rs : C \to W(A,B)$ est lin\Ž aire. Plus
pr\Ž sis\Ž ment, on a le r\Ž sultat suivant.

\

\su{11.3. Th\Ž or\ me} $<17>$.

\

\sl Soit $A \to C \overset{\rp} \to \to B$ une extension de
l'alg\ bre de Lie $B$ par l'alg\ bre de Lie $A$. Pour
toute section lin\Ž aire $\rs : C \to B$ de $\rp$, la repr\Ž
sentation associ\Ž e $f_\rs : C \to W(A,B)$ est un
homomorphisme \bf injectif \sl de l'alg\ bre de Lie $C$ dans
le produit d'entrelacement $W(A,B)$\rm.

\

\su{D\Ž monstration}

\

\su{1. L'application lin\Ž aire $f_\rs$ est injective}

\

\noi \bf En effet\rm, si $f_\rs(c) = 0$ alors, d'une part,
on a
$\rp c = 0$ et, d'autre part, on a $h_{c,0} \equiv
0$. Or, $h_{c,0} = c - \rs \rp c$ donc $c =
0$.\qed

\

\su{2. L'application $f_\rs$ est un homomorphisme} Cette
longue d\Ž monstration s'\Ž tend sur les pages 48 \ˆ \ 52.

\

Soient $a \in C$ et $b \in C$. Il s'agit de montrer que l'on
a
$$(h_{[a,b]},\rp [a,b]) = ([h_a,h_b] + \rp a\star h_b -
\rp b\star h_a, [\rp a,\rp b])$$
et, puisque $\rp [a,b] = [\rp a,\rp b]$, cela revient \ˆ \
montrer que l'on a
$$h_{[a,b]} = [h_a,h_b] + \rp a\star h_b -
\rp b\star h_a.$$
Afin de rendre les calculs plus faciles, on
introduit les polyn\™ mes homog\ nes  suivants qui
appartiennent \ˆ \
$C[[B]]$ :
$$u_{c,m} = \frac{1}{m!} R^m c \ , \ v_{c,m} =
\sum_{n+r = m} \frac{t_r}{(n+1)!} R^n e R^r c$$ 
puis les s\Ž ries formelles
$$u_c = \sum_m u_{c,m} \ , \ v_c = \sum_m v_{c,m} , \
\text{de sorte que} \ \ h_c = u_c - v_c.$$
On a ainsi
$$h_{[a,b]} = u_{[a,b]} - v_{[a,b]},$$
$$\rp a\star h_b -
\rp b\star h_a = \rp a\star u_b - \rp b\star u_a -
\rp a\star v_b + \rp b\star v_a.$$

\

\su{Rappel des donn\Ž es du calcul} On a
$$\rd_{\rp a,n}(y) = t_n (\rad \ y)^n (\rp a) = t_n D^n \rp
a = t_n \rp R^n a.$$
L'alg\ bre de Lie $C$ agit sur $B$ au travers de
l'application compos\Ž e 
$$C \overset{\rp} \to \to B \overset{\rd} \to \to S(B).$$
\bf On va \Ž tablir, successivement, les identit\Ž s
suivantes\rm
$$u_{[a,b]} = [u_a,u_b].\tag 1$$
$$\rp a\star u_b = [v_a,u_b] \ \ \text{et donc} \ \
\rp b\star u_a = [v_b,u_a].\tag 2$$
$$\rp a\star v_b - \rp b\star v_a = [v_a,v_b] + v_{[a,b]}
.\tag 3$$
\bf D'o\ \ d\Ž coulera \rm ceci :
$$[h_a,h_b] + \rp a\star h_b -
\rp b\star h_a =$$
$$[u_a,u_b] - [u_a,v_b] -
[v_a,u_b] + [v_a,v_b] + [v_a,u_b] - [v_b,u_a]  -
[v_a,v_b] - v_{[a,b]} =
u_{[a,b]} - v_{[a,b]} = h_{[a,b]}.$$
Ce qu'il fallait d\Ž montrer.

\

\su{Voici le d\Ž tail des calculs} \bf Tous les indices qui
interviennent (\rm muets ou non) \bf sont des
entiers $\geq 0$\rm.

\

\noi \bf (1) On a \rm
$$u_{[a,b],m}(y) = \frac{1}{m!}
R^m [a,b].$$
\bf D'autre part, on a \rm
$$[u_a,u_b]_m(y) = \sum_{n+r=m} [u_{a,n},u_{b,r}](y) =
\sum_{n+r=m} \frac{1}{n! r!} [R^na,R^rb].$$ 
La formule de Leibniz donne l'\Ž galit\Ž \ (1) souhait\Ž e.

\

\noi \bf (2) On a \rm
$$(\rp a\star u_b)_m(y) = \sum \Sb j+l = m \\ k \leq l
\endSb \frac{t_j}{(l+1)!} R^k [eR^ja,R^{l-k} b].$$
Le lemme \bf 9.8\rm(1), appliqu\Ž \ \ˆ \ la d\Ž rivation $R$
de l'alg\ bre $C$, donne
$$\sum_{k \leq l}  R^k [eR^j a,R^{l-k}b] = \sum_{i \leq l} 
\binom {l+1}{i+1} [R^i e R^j a, R^{l-i} b],$$
donc
$$\frac{t_j}{(l+1)!}\sum_{k \leq l}  R^k [eR^j a,R^{l-k}b] =
\sum_{r \leq l} \frac{t_j}{(i+1)! r!} [R^i e R^j a,R^r b].$$
\bf D'autre part, on a\rm 
$$[v_a,u_b]_m(y) = \sum_{n+r=m} [v_{a,n},u_{b,r}](y) =
\sum_{i+j+r=m} \frac{t_j}{(i+1)! r!} [R^i e R^j a,R^rb].$$
D'o\ \ d\Ž coule l'\Ž galit\Ž \ (2) annonc\Ž e.

\

\noi \bf (3) Calcul de $(\rp a\star v_b
- \rp b\star v_a)_m(y)$\rm. Pour ce calcul, le plus
compliqu\Ž , on introduit
$$w_{c,n,r}(y) = t_r R^n e R^r c$$
afin d'\Ž crire
$$v_{c,m} = \sum_{n+r = m} \frac{1}{(n+1)!} w_{c,n,r}.$$
\bf On a\rm
$$(\rp R^n a).(R^j e R^r b) = \sum_{k \leq j-1}
R^k [e R^n a,R^{j-k-1} e R^r b] + \sum_{k \leq r-1}
R^j e R^k [e R^n a,R^{r-k-1} b].$$
En utilisant la \sl fermeture \Ž clair\rm, il vient
$$(\rp R^n a).(R^j e R^r b) = \sum_{k \leq j-1}
R^k [e R^n a,R^{j-k-1} e R^r b] + \sum_{k \leq r-1}
R^j \rs  D^k [D^n \rp a,D^{r-k-1} \rp b].$$
En \Ž changeant $a$ et $b$ d'une part, et $n$ et $r$
d'autre part, il vient
$$(\rp R^r b).(R^j e R^n a) = \sum_{k \leq j-1}
R^k [e R^r b,R^{j-k-1} e R^n a] + \sum_{k \leq n-1}
R^j \rs D^k [D^r \rp b,D^{n-k-1} \rp a].$$
\bf D'o\ \rm
$$(\rd_{\rp a,n}.w_{b,j,r} - \rd_{\rp
b,r}.w_{a,j,n})(y) = t_n t_r ((\rp R^n a).(R^j e R^r b) -
(\rp R^r b).(R^j e R^n a)) = $$
$$t_n t_r \sum_{k \leq j-1} (R^k [e R^n
a,R^{j-k-1} e R^r b] -  R^k [e
R^r b,R^{j-k-1} e R^n a])  + $$
$$t_n t_r R^j \rs \left(\sum_{k \leq r-1} D^k
[D^n \rp a,D^{r-k-1} \rp b] - \sum_{k \leq n -1}
D^k [D^r \rp b,D^{n-k-1} \rp a]\right).$$
\bf Ainsi \rm
$$(\rp a\star v_b)_m = \sum_{n+s=m+1}(\rd_{\rp a,n}.v_{b,s})
=
\sum_{n+s=m+1}
\rd_{\rp a,n}
\sum_{j+r=s}
\frac{1}{(j+1)!} w_{b,j,r} =$$
$$\sum_{n+j+r=m+1} \frac{1}{(j+1)!} \rd_{\rp
a,n}.w_{b,j,r}.$$
\bf Autrement dit\rm,
$$(\rp a\star v_b)_m(y) = \sum_{n+j+r = m+1} \frac{t_n
t_r}{(j+1)!} (\rp R^n a).(R^j e R^r b).$$
\bf Donc\rm
$$(\rp a\star v_b - \rp b\star v_a)_m(y) = 
\sum_{n+j+r=m+1}
\frac{t_n t_r}{(j+1)!} ((\rp R^n a).(R^j e R^r b) - (\rp R^r
b).(R^j e R^n a)) =$$
$$\sum_{n+j+r=m+1} \frac{t_n
t_r}{(j+1)!} \sum_{k \leq j-1} (R^k [e R^n
a,R^{j-k-1} e R^r b] -  R^k [e
R^r b,R^{j-k-1} e R^n a]) +$$
$$\sum_{n+j+r=m+1} \frac{t_n
t_r }{(j+1)!} R^j \rs\left(\sum_{k \leq r-1} (D^k
[D^n \rp a,D^{r-k-1} \rp b] - \sum_{k \leq n -1} 
D^k [D^r \rp b,D^{n-k-1} \rp a]
\right).$$
\bf On a \rm
$$[v_a,v_b]_m(y) = \sum_{n+r+i+k = m}
\frac{t_nt_r}{(i+1)!(k+1)!} [R^i e R^n a,R^keR^rb].$$

\noi \bf La diff\Ž rence
$(\rp a\star v_b - \rp b\star v_a - [v_a,v_b])_m(y)$\rm.

\

Pour $m,n,r,$ fix\Ž s, le coefficient de $t_n t_r$ dans
cette diff\Ž rence est
$$\sum_{j = m - n - r + 1} \frac{1}{(j+1)!} \sum_{k
\leq j-1} (R^k [e R^n a,R^{j-k-1} e R^r b] -  R^k [e R^r
b,R^{j-k-1} e R^n a]) +$$
$$\sum_{j= m - n - r +1} \frac{1}{(j+1)!} \sum_{k
\leq r-1} (R^j e R^k [e R^n a,R^{r-k-1} b] -  R^j e R^k
[e R^r b,R^{n-k-1} a])-$$
$$\sum_{n+r+i+k = m}
\frac{1}{(i+1)!(k+1)!} [R^i e R^n a,R^keR^rb].$$
\bf La diff\Ž rence $$\sum_{j = m - n - r + 1}
\frac{1}{(j+1)!}
\sum_{0 \leq k
\leq j-1} R^k [e R^n a,R^{j-k-1} e R^r b] -  R^k [e R^r
b,R^{j-k-1} e R^n a]) - $$
$$\sum_{n+r+i+k = m}
\frac{1}{(i+1)!(k+1)!} [R^i e R^n a,R^keR^rb],$$
est nulle\rm.

\

\noi \bf En effet\rm, 
le premier membre s'\Ž crit
$$\frac{1}{(j+1)!}
\sum_{k \leq j-1} (R^k [e R^n a,R^{j-k-1} e R^r b] + 
R^k [R^{j-k-1} e R^n a,e R^r b]),$$
o\ \ $j = m-n-r+1$.

\

\noi Dans l'identit\Ž \ \bf 9.8\rm(2), en rempla\c cant
$m,a,b,\rD,$ respectivement, par $(j-1),e R^n a,e R^r b,R,$ 
on obtient
$$\frac{1}{(j+1)!}
\sum_{k \leq j-1} (R^k [e R^n a,R^{j-k-1} e R^r b] + 
R^k [R^{j-k-1} e R^n a,e R^r b]) = $$
$$\frac{1}{(j+1)!} \sum_{i \leq j-1} \binom {j+1}{i+1} [R^i
e R^n a,R^{j-i-1} e R^r b] = $$
$$\sum_{i \leq j-1} \frac{1}{(i+1)! (j-i)!} [R^i e R^n
a,R^{j-i-1} e R^r b].$$
Puisque $j - 1 = m - n - r$, le second terme s'\Ž crit
$$\sum_{k = j - i -1}
\frac{1}{(i+1)!(k+1)!} [R^i e R^n a,R^keR^rb]$$
et, en rempla\c cant $k$ par $j-i-1$, 
$$\sum_{i \leq j-1} \frac{1}{(i+1)! (j-i)!} [R^i e R^n
a,R^{j-i-1} e R^r b].\qed$$
\bf On a\rm
$$v_{[a,b],m}(y) = \sum_{j} \frac{t_{m-j}}{(j+1)!} R^j e
R^{m-j} [a,b] = \sum_{j} \frac{t_{m-j}}{(j+1)!} R^j \rs (p
R^{m-j} [a,b]) =$$
$$\sum_{j} R^j \rs \frac{t_{m-j}}{(j+1)!}  D^{m-j} [\rp a,\rp
 b].$$

\noi \bf La diff\Ž rence
$$\sum_{j= m - n - r +1} \frac{t_n t_r}{(j+1)!} R^j \rs
\sum_{k
\leq r-1} D^k [D^n \rp a,D^{r-k-1} \rp b] -  \sum_{k
\leq n -1}D^k
[D^r \rp b,D^{n-k-1} \rp a]) -$$
$$\sum_{j} R^j \rs \frac{t_{m-j}}{(j+1)!}  D^{m-j} [\rp a,\rp
 b]$$
est \Ž galement est nulle\rm.

\

\noi \bf En effet\rm,  d'apr\
s l'identit\Ž
\
\bf 9.9\rm, pour $m$ et $j$ fix\Ž s, on a
$$\sum_{n + r = m  -j + 1}{t_n t_r}
\sum_{k
\leq r-1} D^k [D^n \rp a,D^{r-k-1} \rp b] -  \sum_{k
\leq n -1}D^k
[D^r \rp b,D^{n-k-1} \rp a]) = $$
$$t_{m-j} D^{m-j}[\rp a,\rp
b].$$
L'identit\Ž \ (3) en d\Ž coule. Ce qui ach\ ve la d\Ž
monstration du th\Ž or\ me.\qed \qed

\

\su{11.4. Derni\ res remarques}

\

Avec les notations
ci-dessus, dans la formule $f_\rs(c) = (h_c,\rp c)$,  la s\Ž
rie formelle
$h_c$ peut s'\Ž cire
$$h_c = \left(e^R  - \frac{e^R - 1}{R} \  e
\frac{Re^R}{e^R - 1}\right) c =  \left(\roman{Ad} \ z -
\frac{\roman{Ad} \ z - 1}{\rad \ z} \ \rs \rp \ G(\rad \ z)
\right) c,$$
l'exponentielle $e^R = e^{\rad \ z}$ s'\Ž
crivant, traditionnellement, $\roman{Ad} \ z$.

\

Le plongement $f_\rs$ de $C$ dans $W(A,B)$ d\Ž pend de la
section lin\Ž aire $\rs : B \to C$ choisie et il y a une
grande vari\Ž t\Ž \ de choix pour $\rs$, autant que de
suppl\Ž mentaires du sous-espace vectoriel $A$ dans l'espace
vectoriel $C$. Pour tout $a \in A$, on a 
$$h_a = (\roman{Ad} \ z)a = (\roman{Ad} \ \rs y) a,$$
de sorte que 
$$f_\rs(a) = ((\roman{Ad} \ \rs y) a,0)$$
d\Ž pend \Ž galement de $\rs$.

\

Bien entendu, l'alg\ bre de Lie
$W(A,B)$, en tant que produit semi-direct, est une extension
\bf inessentielle \rm de l'alg\ bre de Lie
$B$ par l'alg\ bre de Lie $A[[B]]$.

\

\su{En guise de conclusion} Le cas o\ \ le corps $K$ est
\bf fini \rm n'a pas \Ž t\Ž \ abord\Ž . Il m\Ž rite, sans
doute, d'\ tre examin\Ž \ attentivement. On obtiendrait
ainsi des produits d'entrelacement \bf finis \rm qui
pourraient pr\Ž senter des aspects combinatoires int\Ž
ressants.

\

\

\

\centerline{\sl Le second auteur assume l'enti\ re
responsabilit\Ž
\ de toutes les erreurs}
\centerline{\sl qui seront relev\Ž es dans ce texte
pour avoir
\Ž chapp\Ž \
\ˆ
\ sa vigilence\rm.}

\

\

\

\

\

\

\

\enddocument